\documentclass[12pt]{amsart}
\usepackage{amssymb}
\usepackage{amscd}
\usepackage{amsmath}

\newcommand{\ga}{{\mathfrak a}}
\newcommand{\gb}{{\mathfrak b}}
\newcommand{\gog}{{\mathfrak g}}
\newcommand{\gh}{{\mathfrak h}}
\newcommand{\gl}{{\mathfrak l}}
\newcommand{\gm}{{\mathfrak m}}
\newcommand{\gn}{{\mathfrak n}}

\newcommand{\gs}{{\mathfrak s}}
\newcommand{\bB}{{\bf B}}
\newcommand{\bD}{{\bf D}}
\newcommand{\bG}{{\bf G}}
\newcommand{\bL}{{\bf L}}
\newcommand{\bM}{{\bf M}}

\newcommand{\bP}{{\bf P}}
\newcommand{\bS}{{\bf S}}
\newcommand{\bT}{{\bf T}}

\newcommand{\bR}{{\bf R}}
\newcommand{\al}{\alpha}

\newcommand{\lam}{\lambda}

\newcommand{\cal}{\mathcal} 

\newcommand{\Nscr}{{\cal N}}
\newcommand{\Bscr}{{\cal B}}
\newcommand{\Oscr}{{\cal O}}
\newcommand{\Vscr}{{\cal V}}
\newcommand{\Uscr}{{\cal U}}
\newcommand{\Iscr}{{\cal I}}
\newcommand{\Wscr}{{\cal W}}
\newcommand{\Xscr}{{\cal X}}

\newcommand{\Dscr}{{\cal D}}

\newcommand{\Lie}{{\rm Lie\,}}
\newcommand{\sh}{{\rm sh\,}}
\newcommand{\Rank}{{\rm Rank\,}}

\newcommand{\gr}{{\rm gr\,}}
\newcommand{\Ann}{{\rm Ann\,}}
\newcommand{\QED}{\par \hspace{15cm}$\blacksquare$ \par}
\newcommand{\pr}{^{\prime}}
\newcommand{\prpr}{^{\prime\prime}}
\newcommand{\st}{\subset}
\newcommand{\Co}{{\mathbb C}} 
\newcommand{\Na}{{\mathbb N}}  

\newcommand{\Pf}{\noindent{\bf Proof.}\par\noindent}
\newcommand{\sr}{\scriptscriptstyle}
\newcommand{\parno}{\par\noindent}
\newcommand{\vb}{\vrule height 14pt depth 7pt}
\newcommand{\ts}{\tabskip 4pt}
\newcommand{\vsa}{\noalign{\vskip-7pt}}
\newcommand{\ssa}{\noalign{\vskip -1pt}}

\newcommand{\ov}{\overline}
\newcommand{\ao}{\stackrel {\rm A}{\leq}}     

\newcommand{\go}{\stackrel {\rm G}{\leq}}     

\newcommand{\geg}{\stackrel {\rm G}{\geq}}
\newcommand{\uar}{\uparrow}
\newcommand{\dar}{\downarrow}
\newcommand{\rar}{\rightarrow}

\newtheorem*{theorem}{Theorem}
\newtheorem*{lemma}{Lemma}
\newtheorem{spec}{Claim}
\newtheorem{specia}{Claim}

\newtheorem*{cor}{Corollary}
\newtheorem*{prop}{Proposition}

\begin{document}
\title[$\bB$-orbits of order 2]
{\bf Description of $\bB-$orbit closures of order 2 in
upper-triangular matrices}
\author{Anna Melnikov}\thanks{Supported in part by the Minerva Foundation, Germany,
Grant No. 8466}
\address{Department of Mathematics,
University of Haifa, Haifa 31905, Israel}
\email{melnikov@math.haifa.ac.il}
\begin{abstract}
Let $\gn_n(\Co)$ be the algebra of strictly upper-triangular
$n\times n$ matrices and $\Xscr_2=\{u\in \gn_n(\Co)\ :\ u^{\sr
2}=0\}$ the subset of matrices of nilpotent order 2. Let
$\bB_n(\Co)$ be the group of invertible upper-triangular matrices
acting on $\gn_n$ by conjugation. Let $\Bscr_u$ be the orbit of
$u\in\Xscr_2$ with respect to this action. Let $\bS_n^2$ be the
subset of involutions in the symmetric group $\bS_n.$ We define a
new partial order on $\bS_n^2$ which gives the combinatorial
description of the closure of $\Bscr_u.$ We also construct an ideal
$\Iscr(\Bscr_u)\st S(\gn^*)$ whose variety $\Vscr(\Iscr(\Bscr_u))$
equals $\ov\Bscr_u.$
\par
We apply these results to orbital varieties of nilpotent order 2 in
$\gs\gl_n(\Co)$ in order to give a complete combinatorial
description of the closure of such an orbital variety in terms of
Young tableaux. We also construct the ideal of definition  of such
an orbital variety up to taking the radical.
\end{abstract}
\maketitle
\section{\bf Introduction}
\subsection{}\label{1.1}
Let $\bG=\bS\bL_n(\Co)$ be the special linear group of degree $n$
over field $\Co$. Let $\bB=\bB_n(\Co)$ be the standard Borel
subgroup of $\bG,$ that is the subgroup of all invertible
upper-triangular matrices.
\par
Take $\gog=\gs\gl_n$ and let $\gog=\gn\oplus\gh\oplus\gn^-$ be its
standard triangular decomposition where $\gn=\gn_n$ is a subalgebra
of strictly upper-triangular matrices, $\gh=\gh_n$ is a subalgebra
of diagonal matrices of trace zero and $\gn^-=\gn_n^-$ is a
subalgebra of strictly lower-triangular matrices in $\gs\gl_n$.
\par
We identify $\gn^-$ with $(\gn^+)^*$ through the Killing form on
$\gog$. For any Lie algebra\ $\ga$, let $S(\ga)$ denote its
symmetric algebra.
\par
For any $u\in\gn$ let $\Oscr_u:=\{AuA^{-1}\ |\ A\in\bG\}$ be its
(adjoint) $\bG$-orbit and $\Bscr_u:=\{BuB^{-1}\ |\ B\in \bB\}$ to be
its (adjoint) $\bB$-orbit.
\par
For any $u\in\gn,$ its only characteristic value is 0. So, $\Oscr_u$
is a nilpotent orbit and every nilpotent orbit is obtained in such a
way. Also, the Jordan form of a nilpotent orbit is described by a
partition of $n.$ Obviously, $\Bscr_u\st\Oscr_u\cap \gn.$ Let
$\Oscr$ be a nilpotent orbit and $v\in\gn.$ We say that $\Bscr_v$ is
associated to $\Oscr$ if $\Bscr_v\st \Oscr.$ It is obvious that
$\Bscr_v$ is associated to $\Oscr$ if and only if $\Oscr_v=\Oscr.$
As shown in ~\cite{B} for $n\geq 6,$ the number of $\Bscr_v$ orbits
associated to $\Oscr$ can be infinite. However, the $\bG$ orbits of
nilpotent order 2 decompose into finitely many $\bB$ orbits and so
are a nice exception.
\par
Let $\Xscr_2$ be the set of elements of nilpotent order 2 in $\gn,$
that is $\Xscr_2:=\{u\in\gn\ | \ u^2=0\}$ where $u^2$ is defined by
associative multiplication of matrices. Let $\bS_n$ be the symmetric
group on $n$ elements. Let $\bS_n^2$ be the set of involutions in
$\bS_n,$ that is $\bS_n^2:=\{\sigma\in\bS_n\ |\ \sigma^2=Id\}.$ As
shown in ~\cite{Mx2} there is a natural bijection between
$\bB$-orbits in $\Xscr_2$ and $\bS_n^2.$
\par
In this paper we define a new partial order on $\bS_n^2$ which
completely describes a $\bB$-orbit closure for $\bB$-orbits in
$\Xscr_2$ via this bijection. The combinatorial description of this
order is natural and very simple. Moreover, for a given $\bB$-orbit
$\Bscr\st \Xscr_2$ we construct (via generators) an ideal
$\Iscr(\Bscr)\in S(\gn^-)$ whose variety $\Vscr(\Iscr(\Bscr))$
equals $\ov\Bscr$ and use it to prove our main Theorem \ref{3.5}.
\par
In addition we apply these results to the study of orbital varieties
in $\Xscr_2.$

On the one hand, the ideas of N. Spaltenstein gave me the
inspiration for the order defined here. On the other hand, we apply
the results on $\bB$-orbits to the study of orbital varieties. So we
would like to give a brief description of orbital varieties.
\subsection{}\label{1.1a}
Orbital varieties derive from the works of N. Spaltenstein
~\cite{Sp,Sp2}, and R. Steinberg ~\cite{St1,St2} written during
their studies of the unipotent variety of a complex semisimple group
$\bf G.$ Translation from the unipotent variety of $\bf G$ to the
nilpotent cone of $\mathfrak g =\rm {Lie}(\bf G),$ gave the notion
of an orbital variety. It is defined as follows.
\par
Consider the adjoint action of $\bf G$ on $\mathfrak g.$ Fix some
triangular decomposition $\mathfrak g=\mathfrak n\bigoplus \mathfrak
h\bigoplus\mathfrak n^-.$ A $\bf G$ orbit $\Oscr$ in $\mathfrak g$
is called nilpotent if it consists of nilpotent elements, that is if
${\Oscr}={\bf G}\, x$ for some $x\in\mathfrak n.$ The intersection
$\Oscr\cap\mathfrak n$ is a Lagrangian subvariety of $\Oscr.$ It is
reducible in general. Its irreducible components are called orbital
varieties associated to $\Oscr.$ Orbital varieties play a key role
in the study of primitive ideals in the enveloping algebra
$U(\mathfrak g).$ They also play an important role in Springer's
Weyl group representations.
\par
The first role above can be detailed as follows. Since $\mathfrak g$
is semisimple, we can identify $\mathfrak g^*$ with $\mathfrak g$
through the Killing form. This identification gives an adjoint orbit
a symplectic structure. Then, by ~\cite{Sp2}, ~\cite{St1} and
~\cite{J}, an orbital variety $\Uscr$ associated to a nilpotent
orbit is a Lagrangian subvariety. Following the orbit method one
would like to attach an irreducible representation of the enveloping
algebra $U(\mathfrak g)$ to $\Uscr.$ This should be a simple highest
weight module. Combining the results of A. Joseph and T. A. Springer
one obtains a one to one correspondence between the set of primitive
ideals of $U(\mathfrak g)$ containing the augmentation ideal of its
centre (thus corresponding to integral weights) and the set of
orbital varieties in $\mathfrak g$ corresponding to Lusztig's
special orbits (see for example ~\cite{B-B}). The picture is
especially beautiful for $\mathfrak g=\mathfrak s\mathfrak l_n.$ In
this case, all orbits are special and by ~\cite{M} the associated
variety of a simple highest (integral) weight module is irreducible.
By ~\cite{B-B} and ~\cite{J} in general orbital variety closures are
the irreducible components of an associated variety of a simple
highest weight module. Whilst for $\mathfrak g=\mathfrak s\mathfrak
l_n$ orbital variety closures are themselves associated varieties
and therefore give a natural geometric understanding of the
classification of primitive ideals. This makes their study
especially interesting.
\subsection{}\label{1.2}
There are two descriptions of orbital varieties. The first one,
valid only for $\gs\gl_n$ was given by Spaltenstein. The second one
is a general construction for any semisimple Lie algebra. It was
obtained by Steinberg. We will consider the first one in what
follows. Both descriptions give a very nice combinatorial
characterization of orbital varieties in $\gs\gl_n$ in terms of
Young tableaux. However, both descriptions are not very satisfactory
from a geometric point of view.
\par
The ultimate description of an orbital variety $\Uscr$ would be via
the ideal $I(\Uscr)\st S(\gn^-)$ of definition of $\ov\Uscr,$ that
is a radical ideal $I$ whose variety $\Vscr(I)$ equals $\ov\Uscr.$
In ~\cite[Lecture 7]{J2} the general form that such an ideal should
take is suggested. Let $N$ be the generic matrix of $\gn^-+1d$. Then
$I(\Uscr)$ should be generated by $a_1, a_2,\cdots,a_k$. Here
$a_i=\gr (m_i),$ where $m_i$ are some minors of $N$ and filtration
is by degree. This suggestion derived from an algorithm for
$\Ann_{U(\gn^-)}\Uscr$ based on the Enright functor ~\cite[8.4]{J}
together with calculations in ~\cite{B}.
\par
Of course, the above does not tackle the difficult question as to
which minors $m_i$ to choose. However, in ~\cite{JM} the ideal of
definition of orbital varieties of codimension 1 in a nilradical was
constructed. This is the simplest non-trivial case since such a
variety closure is a complete intersection, thus we had to find only
one non-trivial generating polynomial.
\par
As explained in \ref{1.1}, here we apply the results for $\bB$
orbits in $\Xscr_2$, to orbital varieties associated to nilpotent
orbits of nilpotent order 2. We use the fact shown in
~\cite[4.13]{Msmith} that an orbital variety $\Uscr\st\Xscr_2$
admits a dense $\bB$-orbit $\Bscr_\Uscr.$ Applying our results to
$\Bscr_\Uscr$ we obtain an ideal $\Iscr(\Uscr):=\Iscr(\Bscr_\Uscr)$
whose variety $\Vscr(\Iscr(\Uscr))$ equals $\ov\Uscr.$ Then the
ideal of definition of $\Uscr$ is the radical of $\Iscr(\Uscr).$
\subsection{}\label{1.3}
Another natural problem connected to studying orbital variety
closures is their description as a union of irreducible varieties
each one lying in some nilpotent orbit. In particular, one can ask
whether an orbital variety closure consists of orbital varieties
only or it includes some other $\bB$ stable varieties.
\par
If $\gog\ne\gs\gl_n$ an orbital variety closure is not necessarily
~\cite{Mclos} a union of orbital varieties. However, the
construction only works outside $\gs\gl_n.$ One can conjecture that
in $\gs\gl_n$ an orbital variety closure is a union of orbital
varieties. This conjecture holds for orbital varieties of Richardson
type as it is shown in \cite{Mrich}. As well it is supported by the
computations for all orbital varieties for $n\leq 6.$
\par
In this paper we show in \ref{3.14} that this conjecture holds for
orbital varieties in $\Xscr_2.$
\par
Given orbital varieties $\Uscr,\Wscr\in \gn$ we say that
$\Wscr\geg\Uscr$ if $\Wscr\st \ov\Uscr.$ In ~\cite{Mnew} we give the
combinatorial description of this order for orbital varieties in
$\Xscr_2$ in terms of Young tableaux. Thus, the result in \ref{3.14}
completes the combinatorial description of an orbital variety
closure in $\Xscr_2$ in terms of Young tableaux. Although an orbital
variety in $\Xscr_2$ has a much simpler structure then an orbital
variety in general even this relatively simple case already
demonstrates the complexity of the theory.
\subsection{}\label{1.4}
The ideas we use to study  $\ov\Bscr_u$ are derived from
Spaltenstein's characterization of orbital varieties in $\gs\gl_n$
by Young tableaux. Let us recall it in short.
\par
A nilpotent orbit $\Oscr_u$ is described by a partition giving the
length of blocks in the Jordan form of $u.$ Suppose the Jordan form
of $u$ consists of $k$ non-zero blocks. We order blocks in
non-increasing order and put $\lam^*_i$ to be the length of $i-$th
block. Then $\Oscr_u$ is defined by the partition
$\lam^*=(\lam^*_1,\ldots,\lam_k^*).$ We can represent $\lam^*$
graphically as a Young diagram -- an array of $k$ rows each starting
at the same place on the left, where the $i$-th row contains
$\lam^*_i$ boxes.
\par
 Since in this paper we consider only  nilpotent elements with the Jordan
blocks of length not more than 2, it is more convenient for us to
take the dual partition to $\lam^*$ that is
$\lam=(\lam_1,\ldots,\lam_m)$ where $\lam_i=\#\{\lam^*_j\in\lam^*\
|\ \lam_j\geq i\}.$ In these terms $\Oscr_u$ is described by Young
diagram $D_\lam$ with $m$ columns with $j$-th column containing
$\lam_j$ boxes. We will denote the corresponding orbit by
$\Oscr_\lam.$

Recall that a Young tableau $T$ associated to the Young diagram
$D_\lam$ is obtained by filling  the boxes of $D_\lam$ with numbers
$1,\ldots,n$ so that the numbers increase in rows from left to right
and in columns from top to bottom. Given a Young tableau $T$
associated to $D_\lam$ its shape is defined to be $\lam$ and denoted
as $\sh(T).$ Given $u\in \gn\cap \Oscr_\lam$ its Young diagram is
again defined to be $\lam$ and denoted as $D_n(u)$, or simply
$D(u).$ Now consider canonical projections $\pi_{1,n-i}:\gn_n\rar
\gn_{n-i}$ acting on a matrix by deleting the last $i$ columns and
the last $i$ rows. Set $D_{n-1}(u):=D(\pi_{1,n-1}(u)),
D_{n-2}(u):=D(\pi_{1,n-2}(u)),\ldots, D_1(u):=D(\pi_{1,1}(u))$ and
$\phi(u):=(D_1(u),D_2(u),\ldots,D_n(u)).$ Put $1$ into the unique
box of $D_1(u).$ For any $i\ : 1\leq i\leq n-1$ note that
$D_{i+1}(u)$ differs from $D_i(u)$ by a single "corner" box. Put
$i+1$ into this box. This gives a bijection from the set of the
chains $(D_1(u),D_2(u),\ldots,D_n(u))$ to the set of standard Young
tableaux $T$ of shape $D_n(u).$ In other words, we view a standard
Young tableau as a chain of Young diagrams. For $u\in\gn$ put
$\varphi(u):=T$ if $T$ corresponds under this bijection to
$\phi(u).$ Set $\nu_{\sr T}:=\{u\in \gn\ |\ \varphi(u)=T\}.$
\par
By Spaltenstein ~\cite{Sp} orbital varieties associated to
$\Oscr_\lam$ are parameterized by standard Young tableaux of shape
$\lam$ as follows. Let $\{T_i\}$ be the set of Young tableaux of
shape $\lambda.$ Set $\Uscr_{T_i}:=\ov\nu_{\sr T_i}\cap
\Oscr_\lambda.$ Then  $\{\Uscr_{T_i}\}$ is the set of orbital
varieties associated to $\Oscr_\lambda.$

\subsection{}\label{1.6}
We would like to push Spaltenstein's construction a little bit
further. In a similar way define projections $\pi_{i,n}:\gn\rar
\gn_{\langle i,n\rangle}$ (cf. \ref{2.2}) acting on a matrix by
deleting the first $i-1$ rows and the first $i-1$ columns. Set
$D^i(u):=D(\pi_{n-i+1,n}(u))$ and
$\theta(u):=(D^1(u),\ldots,D^n(u)).$ Put $n$ into the unique box of
$D^1(u).$ Again, $D^{i+1}(u)$ differs from $D^i(u)$ by a single
"corner" box. Insert $n-i+1$ into the left hand corner by inverting
``jeu de taquin'' (cf. ~\cite{Sch}). This gives a second bijection
from the set of chains $(D^1(u),D^2(u),\ldots,D^n(u))$ to the set of
standard Young tableaux $T$ of shape $D^n(u).$ In other words, we
again view a standard Young tableau as a chain of Young diagrams.
For $u\in\gn$ put $\vartheta(u):=T$ if $T$ corresponds under this
bijection to $\theta(u).$ Again, set $\nu^{\sr T}:=\{u\in \gn\ |\
\vartheta(u)=T\}.$ By symmetry about the anti-diagonal it follows
that any orbital variety $\Uscr$ associated to $\Oscr_\lam$ is
obtained as $\Uscr=\ov{\nu^{\sr T}}\cap\Oscr_\lam.$ Moreover, it is
immediate that $\ov{\nu^{\sr T}}\cap\Oscr_\lam=\Uscr_T.$
\par
On the other hand, $\vartheta(u)$ can differ from $\varphi(u).$ For
standard tableaux $T,S$ of shape $\lam$ put $\nu_{\sr T}^{\sr
S}=(\nu_{\sr T}\cap\nu^{\sr S}).$ Then $\nu_{\sr
T}=\coprod\limits_{\sh S=\lam}(\nu_{\sr T}\cap\nu^{\sr S}).$
Obviously, $\nu_{\sr T}^{\sr S}\st \Uscr_T\cap\Uscr_S.$ Thus, $\dim
\nu_{\sr T}^{\sr S}<\dim \Uscr_T$ whenever $S\ne T.$ Hence,
$\nu_{\sr T}^{\sr T}$ is dense in $\Uscr_T.$
\par
For any $i,j:\ 1\leq i\leq j\leq n$ using jeu de taquin to delete
$1,\ldots, i-1$ from a tableau $T$ and deleting boxes containing
$n,n-1,\ldots, j+1$ from $T$ we define $\pi_{i,j}(T)$ -- a tableau
with $j-i+1$ boxes filled in with $i,\ldots,j.$ Set
$D_{i,j}(T):=\sh(\pi_{i,j}(T)).$ In that way we can consider each
Young tableau as a double chain through
$$T\longleftrightarrow\left[ \begin{array}{cccc}
                                      D_{1,1}(T)&D_{1,2}(T)&\ldots&D_{1,n}(T)\cr
                                                &D_{2,2}(T)&\ldots&D_{2,n}(T)\cr
                                                &          &\ddots&\vdots\cr
                                                &          &      &D_{n,n}(T)\cr\end{array}\right]$$
Correspondingly, we can define the projection $\pi_{i,j}:\gn_n\rar
\gn_{j-i+1}$ acting on a matrix  by deleting  the first $i-1$ rows
and the first $i-1$ columns and the last $n-j$ rows and the last
$n-j$ columns. Set
$$\nu\pr_{\sr T}=\{ u\in\gn\ |\ D(\pi_{i,j}(u))=D_{i,j}(T)\}$$
Obviously, $\nu\pr_{\sr T}\st \nu_{\sr T}^{\sr T}.$ Moreover, using
Steinberg construction one gets that $\Uscr_T=\ov{\nu\pr_{\sr
T}}\cap\Oscr_{\lam}.$ This construction provides us the power rank
conditions introduced in ~\cite{v-L}. In turn these conditions
define an ideal $\Iscr_T\st S(\gn^-)$ such that its variety
$\Vscr(\Iscr_T)$ contains $\ov\Uscr_T.$ However, in general for
$n\geq 6$ this ideal is far away from being an ideal of definition
of $\Uscr_T$ since the variety of the ideal is not irreducible and
can include even another orbital variety associated to the same
nilpotent orbit (cf., for example, ~\cite{v-L}).
\subsection{}\label{1.7}
Now we can explain the way we construct $\Iscr(\Bscr)$ for
$\Bscr\st\Xscr_2.$ Note that adjoint $\bB$ action is compatible
projection, that is $\pi_{i,j}(\ov\Bscr_u)=
\ov\Bscr\pr_{\pi_{i,j}(u)}$ where $\Bscr\pr_{u\pr}$ is a
$\bB_{j-i+1}-$orbit of $u\pr\in\gn_{j-i+1}.$ Thus, we can use the
same construction of double chains not only for orbital varieties
but for $\bB$-orbits as well. Note that a double chain defined by
some Young tableau corresponds in general to some set of
$\bB$-orbits, however, a double chain defined by a $\bB$-orbit will
not correspond in general to a standard Young tableau. This double
chain can be described by power rank conditions exactly in the same
way as in the case of an orbital variety. Thus, as in the case of
orbital varieties we can construct an ideal $\Iscr(\Bscr)\st
S(\gn^-)$ as it is explained in detail in \ref{2.5}, \ref{2.5a}.
Again, the variety $\Vscr(\Iscr(\Bscr))$ always contains $\ov\Bscr$
but in general it contains $\bB-$orbits not from $\ov\Bscr$ as well
(cf. \ref{2.5b}).
\par
However, for $\Bscr\st \Xscr^2$ the variety
$\Vscr(\Iscr(\Bscr))=\ov\Bscr.$ This is the main technical result of
the paper.
\subsection{}\label{1.8}
In this paper we do not touch upon two very important problems
connected to our results. The first problem is whether the ideals we
construct are ideals of definition. To show this we have to show
that they are radical ideals. This is a complex question. For some
cases we can answer positively using technique introduced in
~\cite{Ei}, in particular if $\Bscr$ is associated to the minimal
(non-zero) orbit. Indeed, for the nilpotent orbit of the smallest
non-zero dimension, the ideal constructed in \ref{2.5a} for an
orbital variety associated to it, coincides with the ideal of
definiton for such orbital varieties given in \cite{BJ}. These
computations give us a hope that our ideals are ideals of
definition. However, I am far away from complete answer to this
question.

 The second problem is closely connected to the first
one. It is the construction of a strong quantization of orbital
varieties in $\Xscr_2.$ Again, as it is shown in ~\cite{Jmin} the
problem is far from being easy even for the minimal nilpotent orbit.

Another question, that we do not consider here, which however could
be solved using our construction is the description of intersections
of orbital varieties of nilpotent order 2.
\subsection{}\label{1.9}
The body of the paper consists of three sections. In section 2 we
explain all the background essential in the subsequent analysis to
make the paper self-contained. In section 3 we formulate and prove
the results on $\bB$-orbit closures in $\Xscr_2.$ Finally, in
section 4 we explain the results connected to orbital varieties of
nilpotent order 2 and apply the results of section 3 to complete the
combinatorial description of their closures.
\par
In the end one can find the index of notation in which symbols
appearing frequently are given with the subsection where they are
defined. We hope that this will help the reader to find his way
through the paper.
\parno
{\bf Acknowledgments.}\ \ \ This research was done during my stay as
a guest at the Weizmann Institute. I am very grateful to my host A.
Joseph for the invitation and many fruitful discussions through the
various stages of this work.

I would  also like to express my gratitude to the referee. His
numerous remarks helped to improve the presentation and alter some
proofs. His remarks and the remarks of A. Joseph helped me to bring
this paper to a more clear and, hopefully, readable form.

\section
{\bf Geometric and Combinatorial Preliminaries}
\subsection{}\label{2.1}
Let $P(n)$ denote the set of partitions of $n.$ For all
$\lam:=\{\lam_1\geq\lam_2\geq\cdots\geq \lam_k>0\}\in P(n)$ we
define the Young diagram $D_\lam$ to be the array of k {\bf columns}
of boxes each starting at the first row with the $i$-th column
containing $\lam_i$ boxes. Let $\bD_n$ denote the set of all Young
diagrams with $n$ boxes.

 Recall as well that $\lam^*$ denotes the dual partition.
\par
For example, take $\lam=(4,2,2).$ Then $\lam^*=(3,3,1,1)$ and
$$D_\lam=
\vcenter{ \halign{& \hfill#\hfill \tabskip4pt\cr
\multispan{7}{\hrulefill}\cr \ssa \vb & \ &\vb & \ &\vb & \ &
\ts\vb\cr \vsa \multispan{7}{\hrulefill}\cr \ssa \vb & \ &\vb & \ &
\vb &\ & \ts\vb\cr \vsa \multispan{7}{\hrulefill}\cr \ssa \vb & \ &
\ts\vb\cr \vsa \multispan{3}{\hrulefill}\cr \ssa \vb & \ & \ts\vb\cr
\vsa \multispan{3}{\hrulefill}\cr }}$$

\subsection{}\label{2.1a}
It is convenient to replace $\gs\gl_n$ by $\gog\gl_n.$ This makes no
difference to the nilpotent cone $\Nscr$ whilst the adjoint action
(conjugation) by $\bG={\bf GL}_n$ on $\Nscr$ factors to $\bS\bL_n.$
Let $\bB$ be the (Borel) subgroup of upper-triangular matrices in
$\bG.$ Let $\gb=\Lie(\bB)=\gn\oplus\gh$ be the corresponding
subalgebra of $\gog\gl_n.$
\par
Recall the $\Oscr_u$ defined in \ref{1.1} and the $\Oscr_\lam$
defined in \ref{1.4}. For $u\in\gn\cap \Oscr_\lam$ set
$\sh(u):=\lam$ and $\sh(\Oscr_u):=\sh(u).$
\par
In particular, $\Oscr_{(l,k)}$  denotes the nilpotent orbit with $k$
Jordan blocks of length $2$ and $l-k$ Jordan blocks of length 1.
\subsection{}\label{2.2}
Let $e_{i,j}$ be a matrix having 1 in the $ij$-th entry and $0$
elsewhere. Then $\{e_{i,j}\}_{i,j=1}^n$ is a basis of $\gog.$
\par
Let $R\st \gh^*$ denote the set of non-zero roots, $R^+$ the set of
positive roots corresponding to $\gn$ and $\Pi\st R^+$ the resulting
set of simple roots.
\par
Take $i<j$ and let $\al_{i,j}$ be the root which is the weight of
$e_{i,j}.$ Set $\al_{j,i}=-\al_{i,j}.$ We write $\al_{i,i+1}$ simply
as $\al_i.$ Then $\Pi=\{\al_i\}_{i=1}^{n-1}.$ Moreover,
$\al_{i,j}\in R^+\ \Longleftrightarrow\ i<j.$ One has
$$\al_{i,j}=\begin{cases}\sum\limits_{k=i}^{j-1}\al_k & {\rm if}\ i>j\cr
                        -\sum\limits_{k=i}^{j-1}\al_k & {\rm if}\ i<j\cr\end{cases}$$
Let $X_{\al_{i,j}}:=X_{i,j}:=\Co e_{i,j}$ be the root space defined
by $\al_{i,j}\in R.$
\par
For $\al_i\in \Pi,$ let $\bP_{\al_i}$ be the standard parabolic
subgroup of $\bG$ with $\Lie(P_{\al_i})=\gb\oplus
X_{-\al_i}=\gb\oplus X_{i+1,i}.$
\par
For $1\leq i< j\leq n$ let $\bP_{\langle i,j\rangle} $ be the
standard parabolic subgroup of $\bG$ such that
$\bB,\{\bP_{\al_k}\}_{k=i}^{j-1}\st \bP_{\langle i,j\rangle} ,$
equivalently such that
$$\Lie(\bP_{\langle i,j\rangle} )=\gb\oplus\bigoplus\limits_{i\leq s<t\leq j}X_{t,s}$$
Let $\bM_{\langle i,j\rangle} $ be the unipotent radical of
$\bP_{\langle i,j\rangle} $ and $\bL_{\langle i,j\rangle} $ its Levi
factor. Set $\bB_{\langle i,j\rangle} :=\bB\cap \bL_{\langle
i,j\rangle} .$ Set $\gm_{\langle i,j\rangle} :=\Lie(\bM_{\langle
i,j\rangle} )$ and $\gl_{\langle i,j\rangle} :=\Lie(\bL_{\langle
i,j\rangle} ).$ Set $\gn_{\langle i,j\rangle}
:=\bigoplus\limits_{i\leq s<t\leq j} X_{k,s}$ (that is $\gn_{\langle
i,j\rangle} =\gn\cap \gl_{\langle i,j\rangle} )$). We have
decompositions $\bB=\bM_{\langle i,j\rangle} \ltimes\bB_{\langle
i,j\rangle} $ and $\gn=\gn_{\langle i,j\rangle} \oplus\gm_ {\langle
i,j\rangle} .$ They define projections $\pi_{i,j}:\bB\rar
\bB_{\langle i,j\rangle} $ and $\pi_{i,j}:\gn\rar\gn_{\langle
i,j\rangle} $ which we mentioned in \ref{1.6}. For the matrix
$u=(u_{s,t})_{s,t=1}^n$ where $u\in\gn$ or $u\in\bB$ one has
$$\pi_{i,j}(u)
=\left(\begin{array}{cccc}u_{i,i}&u_{i,i+1}&\ldots& u_{i,j}\cr
                              0  &u_{i+1,i+1}&\ldots&u_{i+1,j}\cr
                    \vdots &\vdots&\ddots&\vdots\cr
                      0& 0&\ldots&u_{j,j}\cr\end{array}\right)
$$
\par
Recall \ref{1.6}. For $u\in \gn$ put
$D_{i,j}(u):=\sh(\pi_{i,j}(u)).$ Let $D_u$ be the resulting
``matrix'' of projections shapes, that is
$$(D_u)_{i,j}=\begin{cases} 0, & {\rm if}\ i\geq j;\cr
                      D_{i,j}(u) & {\rm otherwise}\cr\end{cases}$$
Let $\Bscr_u$ denote the $\bB$ orbit of $u.$

For any $u\in\gn$ and $A\in\bB$ one has
$\pi_{i,j}(AuA^{-1})=\pi_{i,j}(A)\pi_{i,j}(u)\pi_{i,j}(A)^{-1}$ so
that $\sh(\pi_{i,j}(AuA^{-1}))=\sh(\pi_{i,j}(u))$ and
$D_{\Bscr_u}:=D_u$ is well defined. By the discussion in \ref{1.6}
one sees that $D_{i,j-1}(u)$ (as well as $D_{i+1,j}$) is obtained
from $D_{i,j}(u)$ just by decreasing one of the entries of $D_{i,j}$
by 1, that is if $D_{i,j}(u)=(\lam_1,\ldots,\lam_k),$ where
$\lam_k>0,$ then there exists $r:\ 1\leq r\leq k$ such that
$D_{i,j-1}(u)=(\lam_1,\ldots,\lam_r-1,\ldots,\lam_k).$
\par
For example,
$$u=\left(\begin{array}{ccccc}0&1&0&0&0\cr
                              0&0&0&1&0\cr
                              0&0&0&0&1\cr
                              0&0&0&0&0\cr
                              0&0&0&0&0\cr\end{array}\right),\quad D_{\Bscr_u}=D_u=
\left(\begin{array}{ccccc}0&(1,1)&(2,1)&(2,1,1)&(2,2,1)\cr
                          0&   0 &  (2)&  (2,1)&  (2,2)\cr
                          0&   0 &  0  &    (2)&  (2,1)\cr
                          0&   0 &  0  &    0  &    (2)\cr
                          0&   0 &  0  &    0  &      0\cr\end{array}\right).$$
\subsection{}\label{2.2a}
Let $W$ be Weyl group of $(\gog,\gh).$ For $\al\in R^+$ let
$s_{\al}\in W$ be the corresponding reflection. For $\al_i\in\Pi$
set $s_i:=s_{\al_i}.$ We identify $W$ with $\bS_n$ by taking $s_i$
to be the elementary permutation interchanging $i,\, i+1.$
\par
 Set $\bS_{\langle i,j\rangle}:=\langle s_{\al_k}\ :\ i\leq k<j\
\rangle$ which we call a parabolic subgroup of $W=\bS_n.$ Set
$D_{\langle i,j\rangle} :=\{w\in \bS_n\ |\ w(\al_k)\in R^+ \ \forall
i\leq k<j\}.$ Then the multiplication map gives a bijection
$\bS_{\langle i,j\rangle} \times D_{\langle i,j\rangle}\simeq
\bS_n.$ Let $\pi_{i,j}:\bS_n\rar \bS_{\langle i,j\rangle} $ be the
projection onto the first factor.
\subsection{}\label{2.3}
Let us define the decomposition of $\Xscr_2$ into $\bB-$orbits.
\par
Recall that $\bS_n^2$ denotes the subset of involutions in $\bS_n.$
Write $\sigma\in \bS_n^2$ as a product of disjoint cycles of length
2. Order entries inside a given cycle to be increasing. Unless
mentioned to the contrary (in section 4) we order the cycles so that
the first entries increase. Thus,
$\sigma=(i_1,j_1)(i_2,j_2)\ldots(i_k,j_k)$ where $i_s<j_s$ for any
$1\leq s\leq k$ and $i_s<i_{s+1}$ for any $1\leq s<k.$
\par
For $\sigma=(i_1,j_1)(i_2,j_2)\ldots(i_k,j_k)$ we construct
$N_\sigma\in\gn$ as follows
$$(N_\sigma)_{i,j}=\begin{cases} 1, & {\rm if}\ i=i_s, j=j_s\ {\rm for\ some}\ 1\leq s\leq k:\cr
                                 0, & {\rm otherwise}.\cr\end{cases}$$
\par
Note that $N_\sigma$ is the upper-triangular part of the permutation
matrix $M_\sigma.$ Set $\Bscr_\sigma:=\Bscr_{N_\sigma}.$ By
\cite[2.2]{Mx2} one has
\begin{prop} $\Xscr_2=\coprod\limits_{\sigma\in\bS_n^2}\Bscr_{\sigma}.$
\end{prop}
This finiteness property is particular for $\Xscr_2.$ For $n\geq 6$
and $\Xscr_k=\{X\in \gn\ |\ X^k=0\}$ where $k\geq 3$ the number of
$\bB$ orbits in $\Xscr_k$ is infinite in general for any infinite
field as it is shown in \cite{B}.
\par
\subsection{}\label{2.3a}
Given $\sigma=(i_1,j_1)(i_2,j_2)\ldots(i_k,j_k)\in \bS_n^2,$ set
$l(\sigma):=k$ and let $\Oscr_\sigma:=\Oscr_{N_\sigma}$  be the
$\bG$ orbit of $N_\sigma.$
\begin{lemma} For $\sigma\in \bS_n^2$ such that $l(\sigma)=k$ one has $\sh(\Oscr_\sigma)=(n-k,k).$
\end{lemma}
\Pf Indeed, $N_{\sigma}\in \Xscr_2$ and its rank is $k.$ Thus, its
Jordan form contains $k$ blocks of length 2 and $n-k$ blocks of
length 1. \QED
\subsection{}\label{2.4}
Given $\sigma=(i_1,j_1)(i_2,j_2)\ldots(i_k,j_k)$ ordered as in
\ref{2.3}. For $s\ : 2\leq s\leq k$ set
$$r_s:=r(i_s,j_s):=\#\{j_p\ :\ p<s,\ j_p<j_s\}+\#\{j_p\ :\ j_p<i_s\}.$$
For example, take $\sigma=(1,6)(3,4)(5,7).$ Then $l(\sigma)=3$ and
$r_2=0,\ r_3=2+1=3.$
\par
By \cite[3.1]{Mx2} one has
\begin{theorem} For $\sigma=(i_1,j_1)(i_2,j_2)\ldots(i_k,j_k)\in \bS_n^2$ one has
$$\dim \Bscr_\sigma=kn+\sum\limits_{s=1}^k(i_s-j_s)-\sum\limits_{s=2}^k r_s.$$
\end{theorem}
\subsection{}\label{2.5}
Given $\lam=(\lam_1,\ldots,\lam_k)$ and $\mu=(\mu_1,\ldots,\mu_l)$
in $P(n)$ we define $\lam\leq \mu$ if for all $s\ :\ 1\leq s\leq
\min(k,l)$ one has
$$\sum\limits_{i=1}^s\lam_i\leq \sum\limits_{i=1}^s\mu_i.$$
By a result of Gerstenhaber (cf. \cite{He}, for example) one has
$$\ov\Oscr_\lam=\coprod\limits_{\mu\geq \lam}\Oscr_\mu.$$
By the Jordan form,  $\sh(\Oscr_u)=(\lam_1,\ldots,\lam_k)$ if and
only if
$$\Rank u=n-\lam_1,\ \Rank u^2=n-(\lam_1+\lam_2),\ldots,\Rank u^{k-1}=\lam_k,\ \Rank u^k=0$$
Moreover, by the above result of Gerstenhaber  $u\pr\in\ov\Oscr_u,$
if and only if
$$\Rank u\pr\leq n-\lam_1=\sum_{s=2}^k\lam_s,\
\Rank {u\pr}^2\leq n-(\lam_1+\lam_2)=\sum_{s=3}^k\lam_s, \ldots,$$
$$\Rank {u\pr}^{k-1}\leq \lam_k,\ \Rank {u\pr}^k=0$$
  We  call
these the power rank conditions. Each power rank condition can be
translated into polynomials in $S(\gog^*)$ as follows. Let us
identify $S(\gog^*)$ with $S(\gog)$ through the coordinate functions
$$x_{j,i}(e_{r,s})=\begin{cases}1, & {\rm if}\ (s,r)=(i,j),\cr 0, & {\rm otherwise.}\cr
\end{cases}\eqno{(*)}$$
Setting $x_{i,j}=e_{j,i}$ identifies $\gog^*$ with $\gog.$ Now,
$\Rank u^i\leq j$ if and only if every minor of order $j+1$ of $u^i$
is equal to $0.$ Let $N$ be the matrix with indeterminate $x_{i,j}$
in the $ij$-th place . Let $\{P^{(i,j)}_t\}_{t=1}^{{n\choose j}^2}$
be the set of all minors of order $j$ of $N^i.$ Then $\ov\Oscr_\lam$
is the zero variety of the ideal of $S(\gog^*)$ generated by
$\left\{P^{(i,j_i)}_t\right\}_{i=1,\ t=1}^{\ k,\  \ {n\choose
j_i}^2}$, where $j_i=1+\sum\limits_{s=i+1}^k\lambda_s.$
\subsection{}\label{2.5a}
Identify $S(\gn^*)$ with $S(\gn^-)$ through \ref{2.5} $(*)$. Take
$u\in\gn$ and put $\Bscr:=\Bscr_u.$ Let us apply \ref{2.5} to
$\Bscr.$
\par
Let $D_\Bscr$ be the matrix defined in \ref{2.2}. From $D_\Bscr$
define an ideal $I(\Bscr)$ (or simply $I$) of $S(\gn^-)$ as follows.
Suppose $(D_\Bscr)_{i,j}=(\lam_1(i,j),\ldots,\lam_k(i,j)).$ For any
$u\in \Bscr$  one has by definition of $D_\Bscr$ that the power rank
conditions on $\pi_{i,j}(u)$ (cf. \cite{v-L}) are as follows:
$$\Rank(\pi_{i,j}(u))=j-i+1-\lam_1(i,j),\ \Rank((\pi_{i,j}(u))^2)=j-i+1-\lam_1(i,j)-\lam_2(i,j),\ \ldots$$
$$\Rank((\pi_{i,j}(u))^{k-1})=\lam_k(i,j),\quad \Rank((\pi_{i,j}(u))^k)=0.\eqno{(*)}$$
In turn, these power rank conditions translate to polynomials in
$S(\gn^-)$ by \ref{2.5}.

Note that if $u\in \Xscr_2$ then the power rank conditions provide
us polynomials of 2 types: quadratic polynomials given by the
conditions $\pi_{i,j}(u^2)=0$ and polynomials of order
$\Rank(\pi_{i,j}(u))+1$ -- these are the minors of order
$\Rank(\pi_{i,j}(u))+1$ of matrix $\pi_{i,j}(u)$ (cf. \ref{3.4} for
the details).

Note that for any $u\in\gn$ one has
\begin{itemize}
\item[(i)] $\pi_{i,j}(u^r)=(\pi_{i,j}(u))^r$ for  any $i,j\ : 1\leq i<j\leq n$ and $r\in \Na.$
\item[(ii)] $\Rank(\pi_{i,j}(u))\leq \Rank(\pi_{i,j-1}(u))+1$ and
$\Rank(\pi_{i,j}(u))\leq \Rank(\pi_{i+1,j}(u))+1$ for any possible
$i,j.$
\end{itemize}

Let $\Iscr_{(i,j)}(\Bscr)$ be the ideal generated by power rank
conditions of $(D_\Bscr)_{i,j}.$

For $r\ :\ 1\leq r\leq n-1$ and  $i\ :\ 1\leq i\leq n-r$ let us call
the set of elements $(i, i+r)$ of a matrix, its $r-$th set. Let
$\Iscr_r(\Bscr)$ be the ideal generated by the power rank conditions
(translated into corresponding minors) of sets $1,\ldots,r,$ of
$D_\Bscr$ (so that
$\Iscr_r(\Bscr)=\Iscr_{r-1}(\Bscr)+\sum_{i=1}^{n-r}\Iscr_{(i,i+r)}(\Bscr)$
and $\Iscr(\Bscr)=\Iscr_{n-1}(\Bscr)$).
\begin{lemma}
Let $\Bscr$ be a $\bB$ orbit in $\gn.$ Let $\Iscr_{r-1}(\Bscr)$ be
the ideal generated by the power rank conditions of sets
$1,\ldots,r-1$ of $D_\Bscr$ Let $\lam=(\lam_1,\lam_2,\ldots)$ be
either the partition $D_{i,i+r-1}$ or $D_{i+1,i+r}.$ For some $l\in
\Na$ we can write $D_{i,i+r}=(\lam_1,\ldots,\lam_l+1,\ldots).$ Then
the polynomials generated by power rank conditions coming from
$D_{i,i+r}$ for powers less than $l$ belong to the ideal
$\Iscr_{r-1}.$ In particular,
\begin{itemize}
\item[(i)] if $D_{i,i+r}=(\lam_1,\ldots,\lam_k,1)$
then $\Iscr_{(i,i+r)}(\Bscr)\st \Iscr_{r-1}(\Bscr);$
\item[(ii)] if $\lam=(\lam_1,\lam_2)$ and
$D_{i,i+r}=(\lam_1,\lam_2+1)$ (that is $l=2$)
            then $\Iscr_{r-1}+\Iscr_{(i,i+r)}=\Iscr_{r-1}+\langle
            \sum_{s=i+1}^{i+r-1}x_{i+r,s}x_{s,i}\rangle$
\end{itemize}
\end{lemma}
\Pf
 Suppose $D_{i,i+r-1}=(\lam_1,\ldots,\lam_k)$ then
$D_{i,i+r}=(\mu_1,\ldots,\mu_k,)$ where $\mu_l=\lam_l+1$ for some
$l\in \Na$ and $\mu_s=\lam_s$ for any $s\ne l.$ The power rank
conditions give us
$$\Rank((\pi_{i,i+r}(u))^p)= \sum_{t=p+1}^k\mu_t=
\begin{cases} \sum_{t=p+1}^k\lam_t+1, & {\rm if}\ p<l;\cr
              \sum_{t=p+1}^k\lam_t, & {\rm if}\ p\geq l;\cr\end{cases}$$
Define $N$ as in \ref{2.5}. Let $\{m_{i,j}^p(N^s)\}$ be the set of
minors of order $p$ of $(\pi_{i,j}(N))^s,$ where $\pi_{i,j}(N)$ is
obtained by deleting the first $i-1$ rows and columns and the last
$n-j$ rows and columns. Since $N$ is upper-triangular one has
$\pi_{i,j}(N^s)=(\pi_{i,j}(N))^s.$ Then the condition
$\Rank((\pi_{i,i+r}(u)^p))=\sum_{t=p+1}^k\mu_t$ gives the subset of
minors $\{m_{i,i+r}^{1+\sum_{t=p+1}^k\mu_t}(N^p)\}.$ If $p<l$ one
has $m_{i,i+r-1}^{1+\sum_{t=p+1}^k\lam_t}(N^p)\in
\Iscr_{r-1}(\Bscr).$ Taking the expansion of
$m_{i,i+r}^{1+\sum_{t=p+1}^k\mu_t}(N^p)$ by its last column we get
minors $\{m_{i,i+r-1}^{1+\sum_{t=p+1}^k\lam_t}(N^p)\}$ in the sum so
that $m_{i,i+r}^{1+\sum_{t=p+1}^k\mu_t}(N^p)\in \Iscr_{r-1}(\Bscr).$
\QED

Let us determine $\Iscr(\Bscr)$ for the example in \ref{2.2}. Recall
that
$$D_{\Bscr_u}=
\left(\begin{array}{ccccc}0&(1,1)&(2,1)&(2,1,1)&(2,2,1)\cr
                          0& 0   &(2)  &  (2,1)&(2,2)\cr
                          0& 0   &  0  &    (2)&(2,1)\cr
                          0& 0   &  0  &  0    &  (2)\cr
                          0&0&0&0&0\cr\end{array}\right).$$
In our notation we get $\Iscr_1(\Bscr)=\langle P_1=x_{3,2},\,
P_2=x_{4,3},P_3=x_{5,4}\rangle.$ As for the entries of the second
set one has $D_{i,i+2}=(2,1)=(D_{i+1,i+2},1)$ for all  $i\ :\ 1\leq
i\leq 3,$ thus, the second set does not add new generating
polynomials, so that $\Iscr_2(\Bscr)=\Iscr_1(\Bscr).$ Again, since
$D_{1,4}=(2,1,1)=(D_{1,3},1)$ then by our lemma we do not get new
generating polynomials from this condition. The only possible new
generating polynomial produced by the third set obtains from
$D_{2,5}=(2,2)$ and it comes from $\Rank(\pi_{2,5}(N^2))=0.$ This
gives us $P_3=x_{5,3}x_{3,2}+x_{5,4}x_{4,2}$ however
$P_3=P_1x_{5,3}+x_{4,2}P_3$ so that it is in $\Iscr_1(\Bscr).$
Finally, $D_{1,5}=(2,2,1)=(D_{2,5},1)$ so that again by our lemma we
do not get new generating polynomials. Thus,
$\Iscr(\Bscr)=\Iscr_1(\Bscr)=\langle x_{3,2}, x_{4,3},
x_{5,4}\rangle.$
\subsection{}\label{2.5b}
Given two $\bB-$orbits $\Bscr,\Bscr\pr$ applying the result of
Gerstenhaber we get $\Iscr(\Bscr)\st \Iscr(\Bscr\pr)$ if and only if
for any $i<j$ one has $(D_\Bscr)_{i,j}\leq (D_{\Bscr\pr})_{i,j}.$
\par
Obviously, $\ov\Bscr\st \Vscr(\Iscr(\Bscr)).$ However, in general
$\Vscr(\Iscr(\Bscr))\supsetneq \ov\Bscr.$ Power rank conditions are
too weak even to separate two $\bB$-orbits of the same dimension
lying in the same nilpotent orbit $\Oscr.$ For example, consider
$$M=\left(\begin{array}{cccccc}
             0&1&0&0&0&0\cr
             0&0&0&0&0&0\cr
             0&0&0&1&0&0\cr
             0&0&0&0&1&0\cr
             0&0&0&0&0&0\cr
             0&0&0&0&0&0\cr\end{array}\right),\qquad {\rm and} \qquad
  N=\left(\begin{array}{cccccc}
             0&1&0&0&0&0\cr
             0&0&0&0&1&1\cr
             0&0&0&0&0&1\cr
             0&0&0&0&1&0\cr
             0&0&0&0&0&0\cr
             0&0&0&0&0&0\cr\end{array}\right)$$
One has $M,N\in\Oscr_{(3,2,1)}$ and $\dim \Bscr_N=\dim \Bscr_M=11$
thus $\ov\Bscr_N\not\st \ov\Bscr_M.$ On the other hand
$$D_M=\left(\begin{array}{cccccc}
0&(1,1)&(2,1)&(2,2)&(2,2,1)&(3,2,1)\cr 0& 0&
(2)&(2,1)&(2,1,1)&(3,1,1)\cr 0& 0& 0&(1,1)&(1,1,1)&(2,1,1)\cr 0& 0&
0& 0& (1,1)& (2,1)\cr 0& 0& 0& 0& 0& (2)\cr 0& 0& 0& 0& 0&
0\cr\end{array}\right),$$
$$D_N=\left(\begin{array}{cccccc}
0&(1,1)&(2,1)&(3,1)&(3,1,1)&(3,2,1)\cr 0& 0& (2)& (3)& (3,1)&
(3,2)\cr 0& 0& 0& (2)& (2,1)& (2,2)\cr 0& 0& 0& 0& (1,1)& (2,1)\cr
0& 0& 0& 0& 0& (2)\cr 0& 0& 0& 0& 0& 0\cr\end{array}\right).$$ So
that $(D_M)_{i,j}\geq (D_N)_{i,j}$ for all possible $i,j$.Thus,
$\Iscr(\Bscr_M)\subsetneq \Iscr(\Bscr_N)$ so that $\ov\Bscr_N\st
\Vscr(\Iscr(\Bscr_N))\st \Vscr(\Iscr(\Bscr_M)).$ Hence,
$\Vscr(\Iscr(\Bscr_M))\supsetneq \ov\Bscr_M.$
\subsection{}\label{2.6}
From the previous subsections we see that we can define two partial
orders on $\bB$ orbits in $\gn.$ The first one is defined by
inclusion of $\bB$ orbit closures and the second one by inverse
inclusion of the ideals generated by power rank conditions defined
by a $\bB$ orbit. We call the first order geometric and the second
order algebraic. From \ref{2.5b} we see that in general the
algebraic order is a proper extension of the geometric order.

Now consider just $\bB$ orbits of nilpotent order 2 in $\gn.$ They
are in bijection with $\bS_n^2$ so both these orders induce partial
orders on $\bS_n^2.$
\par
The geometric order on $\bS_n^2$ is induced by $\bB$ orbit closures
inclusion, that is for $\sigma, \sigma\pr \in \bS_n^2$ we set
$\sigma\pr \go \sigma$ if $\ov\Bscr_{\sigma\pr}\st\ov\Bscr_\sigma.$
\par
The algebraic order on $\bS_n^2$ is induced by inverse inclusions of
$\Iscr(\Bscr_\sigma),$ that is for $\sigma, \sigma\pr \in \bS_n^2$
we set $\sigma\pr \ao \sigma$ if $\Iscr(\Bscr_{\sigma\pr})\supset
\Iscr(\Bscr_\sigma)$ (or equivalently if
$\Vscr(\Iscr(\Bscr_{\sigma\pr}))\st \Vscr(\Iscr(\Bscr_\sigma))$).

We show in what follows that on $\bS_n^2$ the geometric and the
algebraic orders coincide and we have the precise combinatorial
description of this order.
\subsection{}\label{2.10} Since $\pi_{i,j}$ is a projection it is
continuous, one has
\begin{lemma} For $\bB$-orbit $\Bscr$ and any $i,j\ :\ 1\leq i<j\leq n$ one has
$\pi_{i,j}(\ov\Bscr)\subseteq\ov{\pi_{i,j}(\Bscr)}.$ In particular,
if $\Bscr\pr\st\ov\Bscr$ then
$\pi_{i,j}(\Bscr\pr)\st\ov{\pi_{i,j}(\Bscr)}.$
\end{lemma}
\subsection{}\label{2.11}
We need also the following notation for $\sigma\in\bS_n^2.$ We write
$(i,j)\in\sigma$ if it is one of the cycles of $\sigma.$

We write $\sigma=(i_1,j_1)\ldots(i_k,j_k)$ in increasing order as
described in \ref{2.3}. Set $I(\sigma):=\{i_s\}_{s=1}^k$ and
$J(\sigma):=\{j_s\}_{s=1}^k.$
\par
For $1\leq i<j\leq n$ set $S_{i,j}$ to be the set of indices
$S_{i,j}(\sigma):=\{s\ :\ i\leq i_s\leq j \ {\rm and}\ i\leq j_s\leq
j\}.$ For example, take $\sigma=(1,8)(2,5)(3,4)(6,7)\in \bS_8$. Then
$S_{2,6}(\sigma)=\{2,3\}.$
 One has  $\pi_{i,j}(\sigma)=(i_{s_1},j_{s_1})\ldots (i_{s_q},j_{s_q}),$
where $S_{i,j}=\{s_p\}_{p=1}^q.$ In our example
$\pi_{2,6}(\sigma)=(2,5)(3,4).$
\section
{\bf $\bB$-orbit closures in $\Xscr_2$}
\subsection{}\label{3.1}
For a matrix $u\in\gn$ recall the notion $\Bscr_u$ from \ref{1.1}
and $D_u, D_{\Bscr_u}$ from \ref{2.2}. For $i,j\ :\ 1\leq i< j\leq
n$ set $R_{i,j}(u)=\Rank \pi_{i,j}(u)$ and let $R_u$ denote the
corresponding (upper-triangular) matrix of ranks:
$$(R_u)_{i,j}:=\begin{cases} 0, & {\rm if}\ i\geq j;\cr
                    R_{i,j}(u), & {\rm otherwise}.\cr\end{cases}$$
Obviously, for any $v\in\Bscr_u$ one has $R_v=R_u$ so that we may
write $R_{\Bscr_u}:=R_u$ exactly as we have defined $D_{\Bscr_u}$.
\par
Note that $R_u$ can be read off from $D_u.$ Indeed, denote
$D_{i,j}(u):=(\lam^1_{i,j}(u),\lam^2_{i,j}(u),\ldots)$ then
$R_{i,j}(u)=j+1-i-\lam^1_{i,j}(u).$ In general, the information
encoded in $D_u$ is much reacher than the information encoded in
$R_u.$ However, in the case of $u\in \Xscr_2$ these two matrices
give exactly the same information since
$D_{i,j}(u)=(j+1-i-R_{i,j}(u),\, R_{i,j}(u)).$
\par
Adding one row or one column to a matrix can at most increase its
rank by one, so  for any $u\in\gn$ the matrix $R_u$ has the
following properties
\begin{itemize}
\item[(i)] $(R_u)_{i,j}=0$ if $i\geq j$
\item[(ii)] For $i<j$ one has $(R_u)_{i+1,j}\leq (R_u)_{i,j}\leq (R_u)_{i+1,j}+1$
            and $(R_u)_{i,j-1}\leq (R_u)_{i,j}\leq (R_u)_{i,j-1}+1.$
\end{itemize}
By \ref{2.3} for any $u\in \Xscr_2$ there exist a unique
$\sigma\in\bS_n^2$ such that $\Bscr_u=\Bscr_\sigma.$ We get
$R_u=R_{N_\sigma}$ which we denote  by
$R_\sigma:=R_{N_\sigma}=R_{\Bscr_\sigma}.$ By the definition of
$N_\sigma$  in each row and each column there is not more than one
non-zero entry, thus, $(R_{\sigma})_{i,j}$ is simply the number of
ones in the matrix $\pi_{i,j}(N_\sigma).$
\par
For example, take $\sigma=(1,7)(2,3)(4,6)\in \bS_8.$ Then
$$N_\sigma=\left(\begin{array}{cccccccc}
                 0&0&0&0&0&0&1&0\cr
                 0&0&1&0&0&0&0&0\cr
                 0&0&0&0&0&0&0&0\cr
                 0&0&0&0&0&1&0&0\cr
                 0&0&0&0&0&0&0&0\cr
                 0&0&0&0&0&0&0&0\cr
                 0&0&0&0&0&0&0&0\cr
                 0&0&0&0&0&0&0&0\cr\end{array}\right)\quad{\rm and}\quad
R_{\sigma}=\left(\begin{array}{cccccccc}
                 0&0&1&1&1&2&3&3\cr
                 0&0&1&1&1&2&2&2\cr
                 0&0&0&0&0&1&1&1\cr
                 0&0&0&0&0&1&1&1\cr
                 0&0&0&0&0&0&0&0\cr
                 0&0&0&0&0&0&0&0\cr
                 0&0&0&0&0&0&0&0\cr
                 0&0&0&0&0&0&0&0\cr\end{array}\right)$$
In particular, just by the form of $N_\sigma$ one has that
$(N_{\sigma})_{i,j}=1$ implies
$$(R_\sigma)_{i,k}=\begin{cases}(R_\sigma)_{i+1,k}, & {\rm if}\ k<j;\cr
                     (R_\sigma)_{i+1,k}+1, & {\rm if}\ k\geq j;\cr\end{cases}\qquad {\rm and}\qquad
  (R_\sigma)_{k,j}=\begin{cases}(R_\sigma)_{k,j-1}, & {\rm if}\ k>i;\cr
                     (R_\sigma)_{k,j-1}+1, & {\rm if}\ k\geq i.\cr\end{cases}\eqno{(*)}$$
This implies as well that $(N_\sigma)_{k,i}=(N_\sigma)_{j,k}=0$ for
any $k$ so that
$$(R_\sigma)_{k,i}=(R_\sigma)_{k,i-1}\qquad {\rm and}\qquad (R_\sigma)_{j,k}=
(R_\sigma)_{j+1,k}\qquad {\rm for\ any}\ k\ : 1\leq k\leq
n.\eqno{(**)}$$ Moreover, we get
\begin{lemma} Let $\sigma\in \bS_n^2$
then
$(R_\sigma)_{i,j}=(R_\sigma)_{i+1,j}+1=(R_\sigma)_{i,j-1}+1=(R_\sigma)_{i+1,j-1}+1$
if and only if $(i,j)\in\sigma.$
\end{lemma}
\Pf If $(i,j)\in\sigma$ then $(N_\sigma)_{i,j}=1.$ This implies that
$(R_\sigma)_{i,j}=(R_\sigma)_{i+1,j}+1=(R_\sigma)_{i,j-1}+1.$ Also
$(N_\sigma)_{i,j}=1$ implies $(N_\sigma)_{i,k}=(N_\sigma)_{t,j}=0$
for all $k\ne j$ and $t\ne i.$ Hence, in particular,
$(R_\sigma)_{i+1,j}=(R_\sigma)_{i,j-1}=(R_\sigma)_{i+1,j-1}.$
\par
On the other hand, assume that
$(R_\sigma)_{i+1,j}=(R_\sigma)_{i,j-1}=(R_\sigma)_{i+1,j-1}.$ This
means by the form of $N$ that $N_{i,k}=N_{t,j}=0$ for any $k\ :\
1\leq k\leq j-1$ and for any $t\ :\ i+1\leq t\leq n.$ Then
$(R_\sigma)_{i,j}=(R_\sigma)_{i+1,j}+1$ forces $N_{i,j}=1.$
 \QED

 Thus, by $(*), (**)$ and the lemma, any
$R=R_\sigma$ for $\sigma\in \bS_n^2$ has the following additional
property
\begin{itemize}
\item[(iii)] If $R_{i,j}=R_{i+1,j}+1=R_{i,j-1}+1=R_{i+1,j-1}+1$ then
\begin{itemize}
\item[(a)] $R_{i,k}=R_{i+1,k}$ for any $k<j$ and
$R_{i,k}=R_{i+1,k}+1$ for any $k\geq j;$
\item[(b)] $R_{k,j}=R_{k,j-1}$ for any $k>i$ and
$R_{k,j}=R_{k,j-1}+1$ for any $k\leq i;$
\item[(c)] $R_{j,k}=R_{j+1,k}$ and $R_{k,i}=R_{k,i-1}$
for any $k\ :\ 1\leq k\leq n.$
\end{itemize}
\end{itemize}
\subsection{}\label{3.3}
Let us study the structure of $R_\sigma$ more explicitly.
\begin{lemma} Let $R=R_\sigma$ for some $\sigma\in\bS_n^2.$ For any $i,j\ :\ 1\leq i<j\leq n$
only one of the following is possible
\begin{itemize}
\item[(i)]   $R_{i,j}=R_{i+1,j-1}$ and in that case $R_{i+1,j}=R_{i,j-1}=R_{i,j};$
\item[(ii)]  $R_{i,j}=R_{i+1,j-1}+2$ and in that case $R_{i+1,j}=R_{i,j-1}=R_{i+1,j-1}+1;$
\item[(iii)] $R_{i,j}=R_{i+1,j-1}+1$ and $R_{i+1,j}= R_{i,j-1}=R_{i+1,j-1};$
\item[(iv)]  $R_{i,j}=R_{i+1,j-1}+1$ and $R_{i,j-1}=R_{i+1,j-1}+1,$\
$R_{i+1,j}=R_{i,j};$
\item[(v)]   $R_{i,j}=R_{i+1,j-1}+1$ and $R_{i+1,j}=R_{i+1,j-1}+1,$\
$R_{i,j-1}=R_{i,j}.$
\end{itemize}
In each of the above cases $R$ can be presented graphically as
follows:
$${\rm (i)}\ \begin{array}{ccc}\cdots&p&p\cr
                             \cdots&p&p\cr
                                   &\vdots&\vdots\cr\end{array}\quad\quad\quad\quad {\rm (ii)}\
\begin{array}{ccc}\cdots&p+1&p+2\cr
                  \cdots& p& p+1\cr
                        &\vdots&\vdots\end{array}$$
$${\rm(iii)}\
\begin{array}{ccc}\cdots&p&p+1\cr
                  \cdots&p& p\cr
                        &\vdots&\vdots\end{array};\quad {\rm (iv)}\
\begin{array}{ccc}\cdots&p&p+1\cr
                  \cdots&p&p+1\cr
                        &\vdots&\vdots\end{array},\ {\rm (v)}\
\begin{array}{ccc}\cdots&p+1&p+1\cr
                  \cdots&p&p\cr
                        &\vdots&\vdots\end{array}$$
\end{lemma}
\Pf Consider some $N_\sigma.$ The first case is simply equivalent to
$N_{k,j}=N_{i,s}=0$ for any $k\geq i$ and any $s\leq j.$ The second
case is equivalent to the existence of $k>i$ and of $s<j$ such that
$N_{k,j}=N_{i,s}=1.$ Since $R_{i,j}$ is the number of ones in the
submatrix $N_{(i,j)}$ it is obvious that $R_{i,j}=R_{i+1,j-1}+2$ in
that case. The third case occurs if $N_{i,j}=1$ and, thus,
$N_{k,j}=N_{i,s}=0$ for any $k> i$ and any $s< j$ so that
$R_{i+1,j}= R_{i,j-1}=R_{i+1,j-1}.$ The fourth case occurs if
$N_{i,s}=0$ for all $s\leq j$ and $N_{k,j}=1$ for some $k>i.$
Finally, the fifth case occurs if $N_{k,j}=0$ for all $k\geq i$ and
$N_{i,s}=1$ for some $s<j.$ \QED
\subsection{}\label{3.2}
Let $\bR_n^2\st \bM_n$ denote a subset of matrices satisfying
properties (i)-(iii) of \ref{3.1}.
\begin{prop} $\bR_n^2=\{R_{N_\sigma}\}_{\sigma\in\bS_n^2}.$
\end{prop}
\Pf We have seen already that $R_{\sigma}\in \bR_n^2$ for all
$\sigma\in \bS_n^2.$ On the other hand, consider some $R\in \bR_n^2$
and let $N$ be the matrix defined by
$$N_{i,j}=\begin{cases}1, & {\rm if}\ R_{i,j}=R_{i+1,j}+1=R_{i,j-1}+1=R_{i+1,j-1}+1;\  (*)\cr
                       0, & {\rm otherwise}.\cr\end{cases}$$
By property (i) $N_{i,j}=0$ for any $i,j\ :\ i\geq j.$ Note that if
$N_{i,j}=1$ then $N_{i,k}=0$ for any $k\ne j.$ Indeed,
$R_{i,k}=R_{i+1,k}$ for $k<j$ and $R_{i,k-1}=R_{i+1,k-1}+1$ for
$k>j$ by property iii(a) so that conditions $(*)$ are not satisfied
for any $(i,k)$ for which $k\ne j.$ Exactly in the same way we get
by iii(b) that $N_{k,j}=0$ for any $k\ne i.$ Again, by iii(c) no
$N_{j,k}$ and $N_{k,i}$ can satisfy conditions $(*),$ thus,
$N_{j,k}=N_{k,i}=0.$ Hence, $N=N_\sigma$ for $\sigma$ defined by
non-zero entries of $N.$ It remains to show that $R_\sigma=R.$ This
is a straightforward corollary of (i) and (ii) of \ref{3.1}.
 \QED
For example, $\bR_4^2$ contains 10 matrices
$$R_{\sr(1,2)(3,4)}=\left(\begin{array}{cccc}
0&1&1&2\cr 0&0&0&1\cr 0&0&0&1\cr 0&0&0&0\cr\end{array}\right);\
                         R_{\sr(1,4)(2,3)}=\left(\begin{array}{cccc}
0&0&1&2\cr 0&0&1&1\cr 0&0&0&0\cr 0&0&0&0\cr\end{array}\right);\
                         R_{\sr(1,3)(2,4)}=\left(\begin{array}{cccc}
0&0&1&2\cr 0&0&0&1\cr 0&0&0&0\cr 0&0&0&0\cr\end{array}\right);$$
$$\quad
                         R_{(1,2)}=\left(\begin{array}{cccc}
0&1&1&1\cr 0&0&0&0\cr 0&0&0&0\cr 0&0&0&0\cr \end{array}\right);\quad
                     \quad R_{(1,3)}=\left(\begin{array}{cccc}
0&0&1&1\cr 0&0&0&0\cr 0&0&0&0\cr 0&0&0&0\cr\end{array}\right);\quad
                         R_{(1,4)}=\left(\begin{array}{cccc}
0&0&0&1\cr 0&0&0&0\cr 0&0&0&0\cr 0&0&0&0\cr\end{array}\right);$$
$$                     \quad R_{(2,3)}=\left(\begin{array}{cccc}
0&0&1&1\cr 0&0&1&1\cr 0&0&0&0\cr
0&0&0&0\cr\end{array}\right);\quad\quad
                         R_{(2,4)}=\left(\begin{array}{cccc}
0&0&0&1\cr 0&0&0&1\cr 0&0&0&0\cr 0&0&0&0\cr\end{array}\right); $$
$$                         \quad R_{(3,4)}=\left(\begin{array}{cccc}
0&0&0&1\cr 0&0&0&1\cr 0&0&0&1\cr
0&0&0&0\cr\end{array}\right);\quad\quad\
                         R_{Id}=\left(\begin{array}{cccc}
0&0&0&0\cr 0&0&0&0\cr 0&0&0&0\cr 0&0&0&0\cr\end{array}\right).$$
\subsection{}\label{3.4}
Let $I_2\st S(\gn^-)$ be the ideal generated by $I_2=\langle
\sum_{k=j+1}^{i-1}x_{i,k}x_{k,j}\rangle_{j=1,i=j+1}^{j=n-1,i=n}.$
Then $\Vscr(I_2)=\Xscr_2.$
\par
Given $\sigma\in\bS_n^2$ let $\Bscr_\sigma$ be the corresponding
$\bB$ orbit in $\Xscr_2.$
 $D_{\Bscr_\sigma}$ (cf. \ref{2.2}) defines as in
\ref{2.5a} the ideal $\Iscr(\Bscr_\sigma)$ or simply $\Iscr(\sigma)$
of $S(\gn^-).$ which has two types of generating polynomials. The
polynomials of the first type are quadratic (homogeneous)
polynomials generating $I_2$ since $I_2\st \Iscr(\sigma).$ The
generating polynomials of the second type in $\Iscr_\sigma$ are
determined by $\Rank(\pi_{i,j}(N_\sigma))$ that is by
$(R_\sigma)_{i,j}.$ We translate this fact into the set of minors of
$\pi_{i,j}(N)$ of size $(R_\sigma)_{i,j}+1.$ These polynomials are
homogeneous of degree $(R_\sigma)_{i,j}+1.$ Let us denote the ideal
generated by the polynomials of the second type by $I_1(\sigma).$
Put $\Iscr_{R_\sigma}:=\Iscr_\sigma:= I_2+I_1(\sigma).$ Set
$\Vscr_{R_\sigma}:=\Vscr_\sigma: =\Vscr(\Iscr_\sigma).$
\par
Note that $\pi_{i,j}(N_\sigma)$ provides new  polynomials generating
$I_1(\sigma)$ compared to those obtained from
$\pi_{i+1,j}(N_\sigma)$ and $\pi_{i,j-1}(N_\sigma)$ only in the
situation when $(R_\sigma)_{i,j}=(R_\sigma)_{i+1,j-1}.$ Indeed, by
lemma \ref{2.5a}(ii) if $(R_\sigma)_{i,j}=(R_\sigma)_{i+1,j-1}+1$
then all the new generating polynomials obtained from
$(R_\sigma)_{i,j}$ are in $I_2.$
\par
For example, take $\sigma=(1,2)(3,4)\in \bS_5^2.$ One has
$$I_2=\left\langle\begin{array}{c}
x_{\sr 3,2}x_{\sr 2,1};\ x_{\sr 4,3}x_{\sr 3,2},\ x_{\sr 5,4}x_{\sr
4,3};\ x_{\sr 4,3}x_{\sr 3,1}+x_{\sr 4,2}x_{2,1};\cr x_{\sr
5,4}x_{\sr 4,2}+x_{\sr 5,3}x_{\sr 3,2};\ x_{\sr 5,4}x_{\sr
4,1}+x_{\sr 5,3}x_{\sr 3,1}+x_{\sr 5,2}x_{\sr
2,1}\cr\end{array}\right\rangle$$ and
$$R_{(1,2)(3,4)}=\left(\begin{array}{ccccc}0&1&1&2&2\cr
                          0&0&0&1&1\cr
                          0&0&0&1&1\cr
                          0&0&0&0&0\cr
                          0&0&0&0&0\cr\end{array}\right)$$
So that
$$\begin{array}{rcl}
\Iscr_{(1,2)(3,4)}&=&\left\langle I_2; x_{\sr 3,2};\ x_{\sr 5,4};\
\left\vert\begin{array}{cc} x_{\sr 4,2}&x_{\sr 4,3}\cr
                            x_{\sr 5,2}&x_{\sr 5,3}\cr\end{array}\right\vert
                           =x_{\sr 4,2}x_{\sr 5,3}-x_{\sr 4,3}x_{\sr 5,2}\right\rangle\cr
            &=&\left\langle x_{\sr 3,2};\ x_{\sr 5,4}; \ x_{\sr 4,3}x_{\sr 3,1}+x_{\sr 4,2}x_{\sr 2,1};\
x_{\sr 5,3}x_{\sr 3,1}+x_{\sr 5,2}x_{\sr 2,1};\ x_{\sr 4,2}x_{\sr
5,3}-x_{\sr 4,3}x_{\sr 5,2} \right\rangle\cr
\end{array}$$
Moreover, one can easily check that $\Iscr_{(1,2)(3,4)}$ is a prime
ideal.
\par
By the definition of $\Iscr_\sigma$ one can see at once that
$N_\sigma\in\Vscr_\sigma.$ Hence,
\begin{cor} $\ov\Bscr_\sigma\st \Vscr_\sigma.$
\end{cor}
\subsection{}\label{3.5}
Let us define a partial order on $\bR_n^2.$ Given $R,R\pr\in
\bR_n^2$ we put $R\leq R\pr$ if for any $1\leq i,j\leq n$ one has
$R_{i,j}\leq R\pr_{i,j}.$ Let us define the corresponding partial
order on $\bS_n^2.$ For $\sigma,\sigma\pr\in \bS_n^2$ we put
$\sigma\leq \sigma\pr$ if $R_\sigma\leq R_{\sigma\pr}.$
\par
The aim of this part is to prove
\begin{theorem} For any $\sigma\in \bS_n^2$ one has
$$\ov\Bscr_\sigma=\Vscr_\sigma=\coprod\limits_{\sigma\pr\leq\sigma}\Bscr_{\sigma\pr}.$$
\end{theorem}
To do this we study the inclusion relation on
$\{\Iscr_\sigma\}_{\sigma\in\bS_n^2}$ and show that it coincides
with inverse inclusion relation on
$\{\ov\Bscr_\sigma\}_{\sigma\in\bS_n^2}.$ Then the result obtains as
a corollary.
\subsection{}\label{3.6}
All the power rank conditions in the case of $\Xscr_2$ are defined
just by ranks of the minors of the lower-triangular matrix of
indeterminants (and the common conditions coming from $X^2=0$).
Hence $\sigma\leq \sigma\pr$ (or $R_\sigma \leq R_{\sigma\pr}$)
$\Longleftrightarrow$ $\Iscr_\sigma\supset \Iscr_{\sigma\pr}$
$\Longleftrightarrow$ $\Vscr_\sigma\st\Vscr_{\sigma\pr}.$
\begin{lemma} $\Bscr_{\sigma\pr}\st\ov\Bscr_\sigma$ implies
$\Vscr_{\sigma\pr}\st \Vscr_\sigma.$
\end{lemma}
\Pf
The proof is by induction on $n.$ It holds trivially  for $n=2.$
Assume it holds for $n-1.$
\par
Assume $\Bscr_{\sigma\pr}\st\ov\Bscr_\sigma$ for
$\sigma\pr,\sigma\in\bS_n^2.$
\begin{itemize}
\item[(i)] Since $\Bscr_{\sigma\pr}\st\ov\Bscr_\sigma$ implies
$\Oscr_{\sigma\pr}\st\ov\Oscr_\sigma,$ thus, by Gerstenhaber's
result (cf. \ref{2.5}), $\Rank N_{\sigma\pr}\leq \Rank N_\sigma,$
that is, $(R_{\sigma\pr})_{1,n}\leq (R_{\sigma})_{1,n}.$
\item[(ii)] By \ref{2.10} this implies as well that $\pi_{1,n-1}(\Bscr_{\sigma\pr})\st\ov{\pi_{1,n-1}(\Bscr_\sigma)}.$
 By the
induction assumption this provides $R_{\pi_{1,n-1}(\sigma\pr)}\leq
R_{\pi_{1,n-1}(\sigma)}.$
\item[(iii)] Exactly in the same way we get $\pi _{2,n}(\Bscr_{\sigma\pr})\st\ov{\pi_{2,n}(\Bscr_\sigma)}$
and, thus, $R_{\pi_{2,n}(\sigma\pr)}\leq R_{\pi_{2,n}(\sigma)}.$
\end{itemize}
By (ii) and (iii) we get $(R_{\sigma\pr})_{i,j}\leq
(R_\sigma)_{i,j}$ for any $(i,j)\ne (1,n)$ and by (i) one has
$(R_{\sigma\pr})_{1,n}\leq (R_\sigma)_{1,n}.$ Thus,
$R_{\sigma\pr}\leq R_\sigma.$ \QED
\subsection{}\label{3.7}
 To show that $\Vscr_{\sigma\pr}\st \Vscr_\sigma$
implies $\Bscr_{\sigma\pr}\st \ov\Bscr_\sigma$ let us consider the
following two subsets of $\{\sigma\pr\ :\sigma\pr<\sigma\}:$
\begin{itemize}
\item[(i)] $D_1(\sigma)$ is the set of all $\sigma\pr$ such that
\begin{itemize}
\item[(a)] $\sigma\pr<\sigma;$
\item[(b)] $l(\sigma\pr)=l(\sigma);$
\item[(c)] for any $\sigma\prpr$  such
that $\sigma\geq \sigma\prpr\geq \sigma\pr$  one\ has either
$\sigma\prpr=\sigma$ or $\sigma\prpr=\sigma\pr.$
\end{itemize}
\item[(ii)] $D_2(\sigma)$ is the set of all $\sigma\pr$ such that
\begin{itemize}
\item[(a)] $\sigma\pr<\sigma;$
\item[(b)] $l(\sigma\pr)<l(\sigma);$
\item[(c)] for any $\sigma\prpr$ such that
$l(\sigma\prpr)<l(\sigma)$ and $\sigma>\sigma\prpr\geq \sigma\pr$
one has $\sigma\prpr=\sigma\pr.$
\end{itemize}
\end{itemize}

For any $\sigma\pr\in D_1(\sigma)$ one has
$\Oscr_{\sigma\pr}=\Oscr_\sigma.$ Moreover, $D_1(\sigma)$ is the set
of the maximal elements of length $l(\sigma)$ which are smaller than
$\sigma$ with respect to the order defined in \ref{3.5}.

For every $\sigma\pr\in D_2(\sigma)$ one has
$\ov\Oscr_{\sigma\pr}\subsetneq \ov\Oscr_\sigma.$ Moreover
$D_2(\sigma)$ is the set of the maximal elements of length less than
$l(\sigma)$ which are smaller than $\sigma$ with respect to the
order defined in \ref{3.5}.
\par
We begin with the construction of $D_2(\sigma).$ Then we construct
$D_1(\sigma).$ Both constructions are straightforward but the proofs
are rather long.
\par
Recall from \ref{2.3} that we write
$\sigma=(i_1,j_1)(i_2,j_2)\ldots(i_k,j_k)$ where $i_s<j_s$ for any
$s\ : 1\leq s\leq k$ and $i_s<i_{s+1}$ for any $s\ :\ 1\leq s<k.$
Recall notation $(i,j)\in\sigma$, $I(\sigma)$ and $J(\sigma)$ from
\ref{2.11}. Set
$$M(\sigma):=\{s \ :\ s=1\ {\rm or}\  s>1\ {\rm and}\ j_s>\max\{j_1,\ldots j_{s-1}\}\}.$$
The following lemma is straightforward.
\begin{lemma}
If $s\in M(\sigma)$ then for any $(i,j)\in\sigma$ such that
$(i,j)\ne (i_s,j_s)$ one has that either $i>i_s$ or $j<j_s.$ In
particular, fix $m:=\max \{s\in M(\sigma)\}$ then $j_m=\max \{j\in
J(\sigma)\}.$
\end{lemma}
\subsection{}\label{3.7a}
For any $s\ :\ 1\leq s\leq k$ let
$\sigma_s^-:=(i_1,j_1)\ldots(i_{s-1},j_{s-1})(i_{s+1},j_{s+1})\ldots$
be the involution obtained by omitting $(i_s,j_s).$ From the
definition of $\leq$ one has $\sigma_s^-<\sigma.$ Our aim is to show
that $D_2(\sigma)=\{ \sigma_s^-\ |\ s\in M(\sigma)\}.$ We will first
prove the
\begin{lemma}
$D_2(\sigma)\supseteq\{ \sigma_s^-\ |\ s\in M(\sigma)\}.$
\end{lemma}
\Pf Consider $\sigma=(i_1,j_1)\ldots (i_k,j_k)\in \bS_n^2.$
\par
Assume $(i_1,j_1)=(1,n).$ Then
$D_2(\sigma)=\{\sigma_1^-=(i_2,j_2)\ldots(i_k,j_k)\}$ since for any
$(i,j)\ne (1,n)$ one has $(R_\sigma)_{i,j}=(R_{\sigma_1^-})_{i,j}$
and $(R_{\sigma_1^-})_{1,n}= (R_\sigma)_{1,n}-1.$
\par
Now assume $(i_1,j_1)\ne (1,n).$  If $s\in M(\sigma)$ then
$$(R_{\sigma_s^-})_{p,q}=
\begin{cases}(R_{\sigma})_{p,q}-1, & {\rm if}\ p\leq i_s\ {\rm and}\ q\geq j_s;\cr
             (R_{\sigma})_{p,q},   & {\rm otherwise}.\cr\end{cases}$$
Assume that there exists $\sigma\pr$ such that $l(\sigma\pr)<k$ and
$R_\sigma>R_{\sigma\pr}\geq R_{\sigma_s^-}.$ Since
$(R_{\sigma_s^-})_{p,q}=(R_\sigma)_{p,q}$ whenever $p>i_s$ or
$q<j_s$ we obtain
$(R_{\sigma\pr})_{p,q}=(R_{\sigma_s^-})_{p,q}=(R_\sigma)_{p,q}$
whenever $p>i_s$ or $q<j_s.$

In particular, for any $(i,j)\in \sigma_s^-$ we have either $i>i_s$
or $j<j_s$ by lemma \ref{3.7}. Hence,
$(R_{\sigma\pr})_{p,q}=(R_{\sigma_s^-})_{p,q}=(R_{\sigma})_{p,q}$
for $(p,q)\in \{(i+1,j-1),\ (i,j-1),\ (i+1,j),\ (i,j)\}.$ Further,
since $(N_\sigma)_{i,j}=1$ lemma \ref{3.1} gives
$(R_{\sigma})_{i,j}=(R_{\sigma})_{i+1,j-1}+1=
(R_{\sigma})_{i,j-1}+1=(R_{\sigma})_{i+1,j}+1.$ Thus, also
$(R_{\sigma\pr})_{i,j}=(R_{\sigma\pr})_{i+1,j-1}+1=
(R_{\sigma\pr})_{i,j-1}+1=(R_{\sigma\pr})_{i+1,j}+1,$ which by lemma
\ref{3.1} gives $(i,j)\in\sigma\pr$ whenever $(i,j)\in\sigma_s^-.$
Finally, since $l(\sigma\pr)\leq k-1$ we obtain
$\sigma\pr=\sigma_s^-.$ \QED
\subsection{}\label{3.7b}
\begin{prop}
$D_2(\sigma)=\{ \sigma_s^-\ |\ s\in M(\sigma)\}.$
\end{prop}
\Pf
By lemma \ref{3.7a} it remains to show that for any $\sigma\pr$ such
that $\sigma\pr<\sigma$ and $l(\sigma\pr)<l(\sigma)$ there exists
$s\in M(\sigma)$
 such that $ \sigma\pr\leq\sigma_s^-.$

The proof is by induction on $n.$ It is trivially true for
$\sigma\in\bS_2^2.$ Assume it is true for $\sigma\in \bS_{n-1}^2.$
Take $\sigma=(i_1,j_1)\ldots(i_k,j_k)\in \bS_n^2$ and
$\sigma\pr=(i\pr_1,j\pr_1)\ldots (i\pr_{k\pr},j\pr_{k\pr})\in
\bS_n^2$ such that $\sigma\pr<\sigma$ and $k\pr<k.$ In the induction
step we can consider recursion  either to $\pi_{1,n-1}$
 or to $\pi_{2,n}.$ This involves the interplay between all possible
 $(R_{\sigma\pr})_{i,j}$ where
$i=1,2,\ j=n-1,n.$ The argument is broken into the five claims
below. We formulate the assumptions of each claim both in terms of
$J(\sigma\pr),\ J(\sigma)$ (resp. $I(\sigma\pr),\ I(\sigma)$) and in
terms of $(R_{\sigma\pr})_{1,n-1},\ (R_{\sigma})_{1,n-1}$ (resp.
$(R_{\sigma\pr})_{2,n}, \ (R_{\sigma})_{2,n}$) for greater clarity.
The conclusion of each claim is the existence of $s\in M(\sigma)$
such that $\sigma\pr\leq \sigma_s^-.$
\begin{spec}
Assume that $n\not\in J(\sigma\pr),J(\sigma)$ (resp. $1\not\in
I(\sigma\pr),I(\sigma)$). This is equivalent to
$(R_{\sigma\pr})_{1,n-1}=(R_{\sigma\pr})_{1,n}$ and
$(R_{\sigma})_{1,n-1}=(R_{\sigma})_{1,n}$ (resp.
$(R_{\sigma\pr})_{2,n}=(R_{\sigma\pr})_{1,n}$ and
 $(R_{\sigma})_{2,n}=(R_{\sigma})_{1,n}$.) Then there exists $s\in M(\sigma)$ such that
$\sigma\pr\leq \sigma_s^-.$
\end{spec}
\Pf By symmetry about the anti-diagonal it is enough to consider the
case $n\not\in J(\sigma\pr),J(\sigma)$ that is
$(R_{\sigma\pr})_{1,n-1}=(R_{\sigma\pr})_{1,n}=k\pr$ and
$(R_{\sigma})_{1,n-1}=(R_{\sigma})_{1,n}=k.$ In that case
$\pi_{1,n-1}(\sigma\pr)=\sigma\pr$ and $\pi_{1,n-1}(\sigma)=\sigma.$
By the induction hypothesis there exists $s\in M(\sigma)$ such that
$\sigma>\sigma_s^-\geq \sigma\pr.$ Considering $\sigma,\sigma_s^-$
and $\sigma\pr$ as elements of $\bS_n^2$ we get the corresponding
rank matrices by adding an $n-$th row with zero entries and an
$n-$th column with entries equal to those in the $(n-1)$-th column.
Thus, $\sigma\pr\leq \sigma_s^-$ as elements of $\bS_n^2.$ \QED
\begin{spec}
Assume that $n\not\in J(\sigma\pr)$ and $n\in J(\sigma)$ (resp.
$1\not\in I(\sigma\pr)$ and $1\in I(\sigma)$). This is equivalent to
$(R_{\sigma\pr})_{1,n-1}=(R_{\sigma\pr})_{1,n}$ (resp.
 $(R_{\sigma\pr})_{2,n}=(R_{\sigma\pr})_{1,n}$) and
$(R_{\sigma})_{1,n-1}=(R_{\sigma})_{1,n}-1$ (resp.
$(R_{\sigma})_{2,n}=(R_{\sigma})_{1,n}-1$). Then $\sigma\pr\leq
\sigma_m^-$ where $m=\max\{s\in M(\sigma)\}$ (resp.$\sigma\pr\leq
\sigma_1^-$).
\end{spec}
\Pf
  Again, by symmetry about the anti-diagonal it is enough to
consider the case $(R_{\sigma\pr})_{1,n-1}=(R_{\sigma\pr})_{1,n}$
and $(R_{\sigma})_{1,n-1}=(R_{\sigma})_{1,n}-1.$ In that case
$j_m=n.$ Let us show that $\sigma\pr\leq \sigma_m^-.$ Indeed,
$(R_{\sigma_m^-})_{i,j}=(R_{\sigma})_{i,j}$ whenever $j<n$ so that
$(R_{\sigma\pr})_{i,j}\leq (R_{\sigma_m^-})_{i,j}$ for such pairs.
Since $n\not\in J(\sigma\pr)$ one has
$(R_{\sigma\pr})_{i,n}=(R_{\sigma\pr})_{i,n-1}\leq
(R_{\sigma_m^-})_{i,n-1}=(R_{\sigma_m^-})_{i,n}$ for any $i.$ Thus
$\sigma\pr\leq \sigma_m^-.$
 \QED
 We are left with the case $1\in
I(\sigma\pr)$ and $n\in J(\sigma\pr).$ We have two subcases which we
consider in three claims below.
\begin{spec}
Assume that $(1,n)\in\sigma\pr.$ This is equivalent to
$(R_{\sigma\pr})_{2,n-1}=(R_{\sigma\pr})_{1,n-1}=(R_{\sigma\pr})_{2,n}=k\pr-1.$
Then there exists $s\in M(\sigma)$ such that $\sigma\pr<\sigma_s^-.$
\end{spec}
\Pf

If $1\in I(\sigma)$ then $\sigma\pr\leq \sigma_1^-$ since
\begin{itemize}
\item[(i)] $(R_{\sigma_1^-})_{i,j}=
(R_{\sigma})_{i,j}\geq (R_{\sigma\pr})_{i,j}$ whenever $i>1;$
\item[(ii)] $(R_{\sigma\pr})_{1,j}=(R_{\sigma\pr})_{2,j}\leq(R_{\sigma})_{2,j}
=(R_{\sigma_1^-})_{1,j}$ whenever $j<n;$
\item[(iii)] $(R_{\sigma\pr})_{1,n}=k\pr\leq k-1=(R_{\sigma_1^-})_{1,n};$
\end{itemize}

If $1\not\in I(\sigma)$ then $\pi_{2,n}(\sigma)=\sigma$ and
$\pi_{2,n}(\sigma\pr)= (\sigma\pr)_1^-.$ By the induction hypothesis
there exists $s\in M(\sigma)$ such that $\sigma_s^-\geq
(\sigma\pr)_1^-.$ Let us show that $\sigma_s^-\geq \sigma\pr.$
Indeed,
\begin{itemize}
\item[(i)] $(R_{\sigma_s^-})_{i,j}\geq (R_{(\sigma\pr)_1^-})=(R_{\sigma\pr})_{i,j}$
whenever $i>1;$
\item[(ii)] $(R_{\sigma\pr})_{1,j}=(R_{\sigma\pr})_{2,j}\leq(R_{\sigma_s^-})_{2,j}
=(R_{\sigma_s^-})_{1,j}$ whenever $j<n;$
\item[(iii)] $(R_{\sigma\pr})_{1,n}=k\pr\leq k-1=(R_{\sigma_s^-})_{1,n};$
\end{itemize}

Note that in both cases the inequality is strict since in both cases
$(R_{\sigma_s^-})_{2,n}=k-1>k\pr-1=(R_{\sigma\pr})_{2,n}.$
 \QED
  The
remaining case is when $1\in I(\sigma\pr)$ and $n\in J(\sigma\pr),$
but $(1,n)\not\in\sigma\pr.$ We consider two cases according to
whether $n$ is or is not contained in $J(\sigma).$
\begin{spec}
Assume that $1\in I(\sigma\pr)$ and $n\in J(\sigma\pr)$ but
$(1,n)\not\in \sigma\pr.$ This is equivalent to
$(R_{\sigma\pr})_{2,n-1}=k\pr-2,\
(R_{\sigma\pr})_{1,n-1}=(R_{\sigma\pr})_{2,n}=k\pr-1.$ Assume also
that $n\not\in J(\sigma)$. This is equivalent to
$(R_{\sigma})_{1,n-1}=(R_\sigma)_{1,n}$. Then there exists $s\in
M(\sigma)$ such that $\sigma\pr<\sigma_s^-.$
\end{spec}
\Pf If $n\not\in J(\sigma)$ then $\pi_{1,n-1}(\sigma)=\sigma$ and
$l(\pi_{1,n-1}(\sigma))-l(\pi_{1,n-1}(\sigma\pr))=k-(k\pr-1)\geq 2.$
This inequality permits us to use the induction hypothesis twice.
Thus, there exists $s\in M(\sigma)$ and $t\in M(\sigma_s^-)$ such
that $\pi_{1,n-1}(\sigma\pr)\leq (\sigma_s^-)_t^-< \sigma_s^-.$ Let
$i_m$ be $\min\{i_s,i_t\}$ then $m\in M(\sigma)$ and
$(R_{(\sigma_s^-)_t^-})_{i,j}\leq (R_{\sigma_m^-})_{i,j}$ for all
possible $(i,j).$ In particular, $(R_{\sigma\pr})_{i,j}\leq
(R_{\sigma_m^-})_{i,j}$ whenever $j<n.$ Considering $\sigma_m^-,
(\sigma_s^-)_t^-$ as elements of $\bS_n^2$ we get
$$(R_{\sigma_m^-})_{l,n}=\begin{cases}(R_\sigma)_{l,n}, & {\rm if}\ l>i_m;\cr
                      (R_{(\sigma_s^-)_t^-})_{l,n-1}+1, & {\rm if}\ l\leq i_m.\cr
                      \end{cases}$$
On the other hand $(R_{\sigma\pr})_{l,n}\leq
(R_{\sigma\pr})_{l,n-1}+1\leq (R_{(\sigma_s^-)_t^-})_{l,n-1}+1$ for
any $l$ by definition of $(\sigma_s^-)_t^-$ and the construction of
the rank matrices. Thus, for $l\leq i_m$ we get
$(R_{\sigma\pr})_{l,n}\leq
(R_{\sigma_m^-})_{l,n}=(R_{(\sigma_s^-)_t^-})_{l,n-1}+1.$ Finally,
for $l>i_m$ we get $(R_{\sigma\pr})_{l,n}\leq
(R_{\sigma_m^-})_{l,n}=(R_\sigma)_{l,n}.$ Altogether, this gives
$\sigma\pr<\sigma_m^-.$ Note that the inequality is strict since
$(R_{\sigma_m^-})_{1,n-1}=k-1>(R_{\sigma\pr})_{1,n-1}=k\pr-1.$
\QED
\begin{spec}
Assume that $1\in I(\sigma\pr)$ and $n\in J(\sigma\pr)$ but
$(1,n)\not\in \sigma\pr.$ This is equivalent to
$(R_{\sigma\pr})_{2,n-1}=k\pr-2,\
(R_{\sigma\pr})_{1,n-1}=(R_{\sigma\pr})_{2,n}=k\pr-1.$ Assume also
that $n\in J(\sigma).$ This is equivalent to
$(R_{\sigma})_{1,n-1}=(R_\sigma)_{1,n}-1$. Then there exists $s\in
M(\sigma)$ such that $\sigma\pr\leq\sigma_s^-.$
\end{spec}
\Pf
 Again, consider $\pi_{1,n-1}(\sigma),
\pi_{1,n-1}(\sigma\pr).$ Note that
$\pi_{1,n-1}(\sigma)>\pi_{1,n-1}(\sigma\pr).$ We also have
$l(\pi_{1,n-1}(\sigma))=k-1>l(\pi_{1,n-1}(\sigma\pr))=k\pr-1.$ Thus,
by induction hypothesis there exists $s\in M(\pi_{1,n-1}(\sigma))$
such that $\pi_{1,n-1}(\sigma\pr)\leq (\pi_{1,n-1}(\sigma))_s^-.$
Let $r=\max\{t\in M(\sigma)\}$ then $j_r=n.$ Let
$m=\min\{i_s,i_r\}.$ Again, $m\in M(\sigma)$ by lemma \ref{3.7} and
$(R_{\sigma\pr})_{i,j}\leq (R_{\sigma_m^-})$ whenever $j<n.$ Note
also
$$(R_{\sigma_m^-})_{l,n}=\begin{cases}(R_\sigma)_{l,n}, & {\rm if}\ l>i_m;\cr
                            (R_{\sigma_s^-})_{l,n-1}+1, & {\rm if}\ l\leq i_m.\cr\end{cases}$$
Thus, for $l>i_m$ one has $(R_{\sigma\pr})_{l,n}\leq
(R_\sigma)_{l,n}=(R_{\sigma_m^-})_{l,n}$ and for $l\leq i_m$ one has
$(R_\sigma)_{l,n}\leq (R_\sigma)_{l,n-1}+1\leq
 (R_{\sigma_s^-})_{l,n-1}+1=(R_{\sigma_m^-})_{l,n}.$
Altogether this gives $\sigma\pr\leq\sigma_m^-.$
 \QED
In claims 1-5 we have considered all the possible cases described in
lemma \ref{3.3}. Therefore, the proof is completed.
 \QED
  For example, take
$\sigma=(1,6)(2,3)(4,5)(7,8)(9,12)(10,11).$ Then
$$D_2(\sigma)=\left\{\begin{array}{c}(2,3)(4,5)(7,8)(9,12)(10,11),\ (1,6)(2,3)(4,5)(9,12)(10,11),\cr
(1,6)(2,3)(4,5)(7,8)(10,11)\cr\end{array}\right\}$$
\subsection{}\label{3.8b}
Now we construct $D_1(\sigma).$ It consists of four types of
elements. In the next four subsections  we define these types and
show that an element of the given type is in $D_1(\sigma)$. Then in
subsection \ref{3.8a} we show that $D_1(\sigma)$ consists of these
four types of elements only.

In the next five sections we set $\sigma=(i_1,j_1)\ldots (i_k,j_k)$
where $i_s<i_{s+1}$ for any $s:\ 1\leq s<k$ and $i_s<j_s$ for any
$s:\ 1\leq s\leq k.$
 \begin{lemma}. Let $\sigma,\sigma\pr,\sigma\prpr\in \bS_n^2$ be
 such that $\sigma\geq\sigma\pr\geq \sigma\prpr$ and
 $l(\sigma)=l(\sigma\pr)=l(\sigma\prpr).$
 \begin{itemize}
 \item[(i)] If $I(\sigma)=I(\sigma\prpr)$ then
 $I(\sigma\pr)=I(\sigma).$
 \item[(ii)] If $J(\sigma)=J(\sigma\prpr)$ then
 $J(\sigma\pr)=J(\sigma).$
 \end{itemize}
 \end{lemma}
 \Pf
 (i)\ \
 Assume that $I(\sigma)=I(\sigma\prpr)$. Then by definition of a
 rank  matrix $(R_\sigma)_{i,n}=(R_{\sigma\prpr})_{i,n}$ for all
 $i.$ Since $(R_\sigma)_{i,n}\leq (R_{\sigma\pr})_{i,n}\leq (R_{\sigma\prpr})_{i,n}$
 this forces $(R_\sigma)_{i,n}=(R_{\sigma\pr})_{i,n}$ for all $i$.
 Thus $(R_{\sigma\pr})_{i,n}=(R_{\sigma\pr})_{i+1,n}+1$
 $\Longleftrightarrow$
 $(R_{\sigma})_{i,n}=(R_{\sigma})_{i+1,n}+1.$ Hence, $i\in
 I(\sigma)$ iff $i\in I(\sigma\pr).$

 (ii) is obtained exactly in the same way.

 \QED

\subsection{}\label{3.8}
Let $(i_s,j_s)\in \sigma$. Suppose there exists $m<i_s$ satisfying
both
\begin{itemize}
\item[(i)] $m\not \in I(\sigma)\cup J(\sigma);$
\item[(ii)] Either $m=i_s-1$ or for any $t\ :\ m<t<i_s$ one has
$t\in I(\pi_{1,j_s}(\sigma))\cup J(\pi_{1,j_s}(\sigma)).$
\end{itemize}
Then define $\sigma_{i_s\uar}$ to be obtained from $\sigma$ by just
changing $(i_s,j_s)$ to $(m,j_s).$ By the definition such an $m$ is
unique, if it exists. If such an $m$ does not exist put
$\sigma_{i_s\uar}:=\emptyset.$

For example, take $\sigma=(2,5)(3,4)(7,9)(8,10).$ Then
$\sigma_{2\uar}=(1,5)(3,4)(7,9)(8,10),$\
$\sigma_{3\uar}=\emptyset,$\ $\sigma_{7\uar}=(2,5)(3,4)(6,9)(8,10)$
and $\sigma_{8\uar}=(2,5)(3,4)(7,9)(6,10).$

If $\sigma_{i_s\uar}\ne \emptyset$ then
$$
(R_{\sigma_{i_s}})_{i,j}=
\begin{cases}
(R_\sigma)_{i,j}-1, & {\rm if}\ m<i\leq i_s\ {\rm and}\ j\geq
                                                         j_s;\cr
(R_\sigma)_{i,j},    & {\rm otherwise}; \cr
\end{cases}\eqno(*)$$
Thus, $\sigma_{i_s\uar}<\sigma.$
\begin{lemma}
Take $\sigma=(i_1,j_1)\ldots(i_k,j_k).$ Let $s\ :\ 1\leq s\leq k$ be
such that $\sigma_{i_s\uar}\ne\emptyset.$ Then $\sigma_{i_s\uar}\in
D_1(\sigma)$.
\end{lemma}
\Pf
 The proof is by induction on $n.$
 The minimal $n$ which permits such a situation to occur is $n=3.$
In $\bS_3^2$ there  is a unique case of this type: $\sigma=(2,3)$
and $\sigma_{2\uar}=(1,3).$ It is obvious that $\sigma_{2\uar}\in
D_1(\sigma).$ Assume the claim is true for $n-1.$

Suppose $\sigma_{i_s\uar}\ne\emptyset$ for some $s.$ Let
$\sigma\pr=(i\pr_1,j\pr_1)\ldots (i\pr_k,j\pr_k)$ be such that
$\sigma\geq \sigma\pr\geq \sigma_{i_s\uar}.$  We must show that
either $\sigma\pr=\sigma$ or $\sigma\pr=\sigma_{i_s\uar}.$

Note that by definition $J(\sigma)=J(\sigma_{i_s\uar}).$ By lemma
\ref{3.8b} this implies $J(\sigma\pr)=J(\sigma).$ Set
$J:=J(\sigma).$

 Assume first that $n\not\in
J$. Passing to $\pi_{1,n-1}$ we get $\pi_{1,n-1}(\sigma)=\sigma,$\
$\pi_{1,n-1}(\sigma\pr)=\sigma\pr$\ and
$\pi_{1,n-1}(\sigma_{i_s\uar})=\sigma_{i_s\uar}$ so that by the
induction hypothesis either $\sigma\pr=\sigma$ or
$\sigma\pr=\sigma_{i_s\uar}.$

Now assume that $n\in J.$ We have to consider two cases: either
$j_s=n$ or $j_s<n.$
\begin{itemize}
\item[(i)] Assume first that $j_s=n.$ Then
$\pi_{1,n-1}(\sigma_{i_s\uar})=\pi_{1,n-1}(\sigma)$ so that
$(R_\sigma)_{i,j}=(R_{\sigma_{i_s\uar}})_{i,j}=(R_{\sigma\pr})_{i,j}$
whenever $j<n.$ This is equivalent to $(i_r,j_r)\in \sigma$ if and
only if $(i_r,j_r)\in\sigma_{i_s\uar}$ and thus also if and only if
$(i_r,j_r)\in\sigma\pr$ whenever $r\ne s.$  Suppose
$(i\pr_s,n)\in\sigma\pr.$ Note that
$I(\sigma)\setminus\{i_s\}=I(\sigma_{i_s\uar})\setminus\{m\}=
I(\sigma\pr)\setminus \{i_s\pr\}.$ Denote this common set by $\hat
I:=I(\sigma)\setminus\{i_s\}.$

 By $(*)$ and the structure of the rank
matrices
$(R_{\sigma_{i_s\uar}})_{l,n}=(R_{\sigma})_{l,n}=(R_\sigma)_{l,n-1}$
whenever $l>i_s$ and
$(R_{\sigma_{i_s\uar}})_{l,n}=(R_{\sigma})_{l,n}=
(R_\sigma)_{l,n-1}+1$ whenever $l\leq m.$ These equalities force
$(R_{\sigma\pr})_{l,n}=(R_{\sigma\pr})_{l,n-1}$ whenever $l>i_s$ and
$(R_{\sigma\pr})_{l,n}= (R_{\sigma\pr})_{l,n-1}+1$ whenever $l\leq
m.$ Thus, by the structure of the rank matrices we get $m\leq
i\pr_s\leq i_s.$ Further, by the definition of $\sigma_{i_s\uar}$
for any $t\ :\ m<t<i_s$ one has $t\in \hat I\cup J$ so that either
$i_s\pr=i_s$ (that is $\sigma\pr=\sigma$) or $i_s\pr=m$ (that is
$\sigma\pr=\sigma_{i_s\uar}$).

\item[(ii)] Now assume that $j_s<n.$
\begin{itemize}
\item[(a)]
Suppose $(i_r,n)\in\sigma$ (then  $(i_r,n)\in\sigma_{i_s\uar}$ as
well) and $(i\pr_r,n)\in \sigma\pr.$ Let us show that $i\pr_r=i_r.$
Indeed, by definition of $\sigma_{i_s\uar}$ one has either $i_r>i_s$
or $i_r<m.$ In both cases one has by $(*)$ that
$(R_{\sigma_{i_s\uar}})_{i_r,n-1}=(R_\sigma)_{i_r,n-1}$ and
$(R_{\sigma_{i_s\uar}})_{i_r,n}=(R_\sigma)_{i_r,n}=(R_\sigma)_{i_r,n-1}+1.$
One also has
$(R_{\sigma_{i_s\uar}})_{i_r+1,n}=(R_{\sigma_{i_s\uar}})_{i_r,n}-1=
(R_\sigma)_{i_r,n}-1=(R_\sigma)_{i_r+1,n}.$ These equalities force
$(R_{\sigma\pr})_{i_r+1,n}=(R_{\sigma\pr})_{i_r,n-1}=(R_{\sigma\pr})_{i_r,n}-1.$
By the structure of the rank matrices this implies $i\pr_r=i_r.$

\item[(b)] Since $j_s<n$ one has
$\pi_{1,n-1}(\sigma_{i_s\uar})= (\pi_{1,n-1}(\sigma))_{i_s\uar}.$
Thus, by the induction hypothesis either
$\pi_{1,n-1}(\sigma\pr)=\pi_{1,n-1}(\sigma)$ or
$\pi_{1,n-1}(\sigma\pr)=\pi_{1,n-1}(\sigma_{i_s\uar})$. By (a) this
implies that either $\sigma\pr=\sigma$ or
$\sigma\pr=\sigma_{i_s\uar}.$
\end{itemize}
\end{itemize}
\QED

\subsection{}\label{3.8i} Again, let $(i_s,j_s)\in \sigma.$
 Suppose there exists $m>j_s$ satisfying both
\begin{itemize}
\item[(i)] $m\not \in I(\sigma)\cup J(\sigma);$
\item[(ii)] Either $m=j_s+1$ or for any $t\ :\ j_s<t<m$ one has
$t\in I(\pi_{i_s,n}(\sigma))\cup J(\pi_{i_s,n}(\sigma)).$
\end{itemize}
Then define $\sigma_{j_s\rar}$ to be obtained from $\sigma$ just by
changing $(i_s,j_s)$ to $(i_s,m).$ Again, such an $m$ is unique, if
it exists. If such an $m$ does not exist put
$\sigma_{j_s\rar}:=\emptyset.$

For example, take $\sigma=(2,5)(3,4)(7,9)(8,10)\in \bS_{11}^2.$ Then
$\sigma_{5\rar}=(2,6)(3,4)(7,9)(8,10),\ \sigma_{4\rar}=\emptyset,\
\sigma_{9\rar}=(2,5)(3,4)(7,11)(8,10)$ and
$\sigma_{10\rar}=(2,5)(3,4)(7,9)(8,11).$

By symmetry of the rank matrix about the anti-diagonal if
$\sigma_{j_s\rar}\ne\emptyset$ then $\sigma_{j_s\rar}<\sigma$ and
$$
(R_{\sigma_{j_s\rar}})_{i,j}=
\begin{cases}
(R_\sigma)_{i,j}-1, & {\rm if}\ j_s\leq j<m\ {\rm and}\ i\leq
i_s;\cr (R_\sigma)_{i,j},    & {\rm otherwise}; \cr
\end{cases}\eqno(*\pr)
$$
Exactly as in \ref{3.8} we get the
\begin{lemma}
Take $\sigma=(i_1,j_1)\ldots(i_k,j_k).$ Let $s\ :\ 1\leq s\leq k$ be
such that $\sigma_{j_s\rar}\ne\emptyset.$ Then $\sigma_{j_s\rar}\in
D_1$.
\end{lemma}

\subsection{}\label{3.8ii}
 Consider $(i_s,j_s)\in\sigma$ for $s\geq 2$ and suppose there exists
$j_r<i_s$ satisfying one of the conditions: either $ j_r=i_s-1$ or
for any $t\ :\ j_r<t<i_s$ one has $t\in I(\pi_{i_r,j_s}(\sigma))
\cup J(\pi_{i_r,j_s}(\sigma)).$ Then define $\sigma_{j_r,i_s}$ to be
obtained from $\sigma$ just by changing the cycles $(i_r,j_r),\
(i_s,j_s)$ to the cycles $(i_r,i_s),\ (j_r,j_s).$

Note that for a given $i_s$ there can exist a few possible $j_r.$
Let us denote the set of all possible $\sigma_{j_r,i_s}$  for a
given $i_s$ by $C_{i_s\uar\rar}(\sigma).$ If there is no $r<s$
 such that $j_r$ satisfies conditions put $C_{i_s\uar\rar}(\sigma):=\emptyset.$

For example, take $\sigma=(1,3)(2,4)(5,8)(6,7).$ Then
$C_{2\uar\rar}(\sigma)=\emptyset,\ C_{5\uar\rar}(\sigma)=
\{\sigma_{4,5}=(1,3)(2,5)(4,8)(6,7),
\sigma_{3,5}=(1,5)(2,4)(3,8)(6,7)\}$ and $C_{6\uar\rar}=\emptyset.$

When $\sigma_{j_r,i_s}$ is defined one has
$$
(R_{\sigma_{j_r,i_s}})_{i,j}=
\begin{cases}
(R_{\sigma})_{i,j}-1, &
       {\rm if}\ i\leq i_r \ {\rm and}\ j_r\leq j<i_s;\cr
(R_{\sigma})_{i,j}-1, &
       {\rm if}\ i_s\leq i<j_r\ {\rm and}\ j\geq j_s;\cr
(R_{\sigma})_{i,j},   & {\rm otherwise};\cr
\end{cases} \eqno(**)$$
Thus, $\sigma_{j_r,i_s}<\sigma.$
\begin{lemma}
Take $\sigma=(i_1,j_1)\ldots(i_k,j_k).$ Let $s\ :\ 2\leq s\leq k$
and $r\ :\ r<s$ be such that $\sigma_{j_r,i_s}$ exists. Then
$\sigma_{j_r,i_s}\in D_1$.
\end{lemma}
\Pf The proof is by induction on $n.$
 The minimal $n$ which permits such a situation to occur is $n=4.$
In $\bS_4^2$ there  is a unique case of this type:
$\sigma=(1,2)(3,4)$ and $\sigma_{2,3}=(1,3)(2,4).$ It is obvious
that $\sigma_{2,3}\in D_1(\sigma).$ Assume the claim is true for
$n-1$.

Suppose $\sigma_{j_r,i_s}\in C_{i_s\uar\rar}(\sigma)$. Let
$\sigma\pr=(i_1\pr,j_1\pr)\ldots (i\pr_k,j\pr_k)$ be such that
$\sigma_{j_r,i_s}\leq \sigma\pr\leq \sigma.$ We must show either
$\sigma\pr=\sigma$ or $\sigma\pr=\sigma_{j_r,i_s}.$

Note that $I(\sigma)\cup J(\sigma)=I(\sigma_{j_r,i_s})\cup
J(\sigma_{j_r,i_s}).$ In particular, one has $1\in I(\sigma)$ if and
only if $1\in I(\sigma_{j_r,i_s}).$ By $(**)$
$(R_\sigma)_{2,n}=(R_{\sigma_{j_r,i_s}})_{2,n},$ thus,
$(R_{\sigma\pr})_{2,n}=(R_{\sigma})_{2,n}.$ Also,
$(R_\sigma)_{2,n}=k-1$ if and only if $1\in I(\sigma).$ Hence, $1\in
I(\sigma\pr)$ $\Longleftrightarrow$ $1\in I(\sigma).$

Assume first that $1\not\in I(\sigma).$ Passing to $\pi_{2,n}$ we
get $\pi_{2,n}(\sigma)=\sigma,$\ $\pi_{2,n}(\sigma\pr)=\sigma\pr$\
and $\pi_{2,n}(\sigma_{j_r,i_s})=\sigma_{j_r,i_s}$ so that by the
induction  hypothesis either $\sigma\pr=\sigma$ or
$\sigma\pr=\sigma_{j_r,i_s}.$

Now, let $1\in I(\sigma).$
\begin{itemize}
\item[(i)] Assume first that
$i_r>1.$ Then $\pi_{2,n}(\sigma_{j_r,i_s})=
(\pi_{2,n}(\sigma))_{j_r,i_s},$ thus, by the induction hypothesis
$\pi_{2,n}(\sigma_{j_r,i_s})\in D_1((\pi_{2,n}(\sigma)))$ so that
either $\pi_{2,n}(\sigma\pr)=\pi_{2,n}(\sigma)$ or $
\pi_{2,n}(\sigma\pr)=(\pi_{2,n}(\sigma))_{j_r,i_s}.$
\begin{itemize}
\item[(a)] If $\pi_{2,n}(\sigma\pr)=\pi_{2,n}(\sigma)$ then $\sigma\pr<\sigma$
only if $j_1\pr>j_1.$ In that case
$(R_{\sigma\pr})_{1,j_1}=(R_{\sigma})_{1,j_1}-1.$ On the other hand,
since $j_1\not\in [j_r,i_s]$ one has by $(**)$ that
$(R_{\sigma_{j_r,i_s}})_{1,j_1}=(R_\sigma)_{1,j_1},$ which
contradicts the condition $\sigma\pr\geq \sigma_{j_r,i_s}.$ Thus,
$\pi_{2,n}(\sigma\pr)=\pi_{2,n}(\sigma)$ implies $\sigma\pr=\sigma.$
\item[(b)] If $\pi_{2,n}(\sigma\pr)=\pi_{2,n}(\sigma_{j_r,i_s})$ then
$\sigma\pr>\sigma_{j_r,i_s}$ only if $j_1\pr<j_1.$ In that case
$(R_{\sigma\pr})_{1,j_1\pr}=(R_{\sigma_{j_r,i_s}})_{1,j_1\pr}+1.$ On
the other hand since
$\pi_{2,n}(\sigma\pr)=(\pi_{2,n}(\sigma))_{j_r,i_s}$ one has
$j_1\pr\not\in [j_r,i_s]$ so that by $(**)$\
$(R_{\sigma_{j_r,i_s}})_{1,j_1\pr}=(R_\sigma)_{1,j_1\pr}$ which
contradicts the condition $\sigma\pr\leq \sigma.$ Thus,
$\pi_{2,n}(\sigma\pr)=\pi_{2,n}(\sigma_{j_r,i_s})$ implies
$\sigma\pr=\sigma_{j_r,i_s}.$
\end{itemize}
\item[(ii)] Assume that $i_r=1.$ In that case $r=1$ and
$\pi_{2,n}(\sigma_{j_1,i_s})=(\pi_{2,n}(\sigma))_{i_s\uar}.$ Thus,
by Lemma \ref{3.8} either $\pi_{2,n}(\sigma\pr)=\pi_{2,n}(\sigma)$
or $\pi_{2,n}(\sigma\pr)=\pi_{2,n}(\sigma_{j_1,i_s}).$

Assume first that $\pi_{2,n}(\sigma\pr)=\pi_{2,n}(\sigma).$ Then
$\sigma\pr<\sigma$ only if $j\pr_1> j_1.$ By the conditions on the
construction of $\sigma_{j_1,i_s}$
 this implies $j\pr_1>i_s$ so that by $(**)$\
 $(R_{\sigma_{j_1,i_s}})_{1,i_s}=(R_{\sigma\pr})_{1,i_s}+1$
 which contradicts  $\sigma\pr\geq\sigma_{j_1,i_s}.$ Thus,
$\pi_{2,n}(\sigma\pr)=\pi_{2,n}(\sigma)$ implies $\sigma\pr=\sigma.$

 Now assume that
$\pi_{2,n}(\sigma\pr)=\pi_{2,n}(\sigma_{j_1,i_s}).$ Then
$\sigma\pr>\sigma_{j_1,i_s}$ only if  $j\pr_1< i_s.$ Again, by the
conditions on construction of $\sigma_{j_1,i_s}$ this implies
$j\pr_1<j_1$ so that by $(**)$\
$(R_{\sigma\pr})_{1,j\pr_1}=(R_\sigma)_{1,j\pr_1}+1$ which
contradicts $\sigma\pr\leq \sigma.$ Thus,
$\pi_{2,n}(\sigma\pr)=\pi_{2,n}(\sigma_{j_1,i_s})$ implies
$\sigma\pr=\sigma_{j_1,i_s}.$
\end{itemize}
\QED

\subsection{}\label{3.8iii}
In the construction of the last type we use the notion $M(\sigma)$
from \ref{3.7}. For  $(i_s,j_s)\in \sigma$ consider
$\pi_{i_s+1,j_s-1}(\sigma)$ and let
$M_{[i_s,j_s]}(\sigma):=M(\pi_{i_s+1,j_s-1}(\sigma)).$ If
$M_{[i_s,j_s]}(\sigma)\ne\emptyset$ then for  $t\in
M_{[i_s,j_s]}(\sigma)$ define $\sigma_{i_s,i_t}$ to be obtained from
$\sigma$ just by changing the cycles $(i_s,j_s),\ (i_t,j_t)$ to the
cycles $(i_s,j_t),\ (i_t,j_s).$ If $M_{[i_s,j_s]}(\sigma)=\emptyset$
put $C_{i_s\uar\dar}(\sigma):=\emptyset$ otherwise put
$C_{i_s\uar\dar}(\sigma):=\{\sigma_{i_s,i_t}\ :\ t\in
M_{[i_s,j_s]}(\sigma)\}.$

For example, take $\sigma=(1,8)(2,5)(3,4)(6,7).$ Then
$$\begin{array}{lcl}
C_{1\uar\dar}(\sigma)&=& \{\sigma_{1,2}=(1,5)(2,8)(3,4)(6,7),\
\sigma_{1,6}=(1,7)(2,5)(3,4)(6,8)\},\cr
C_{2\uar\dar}(\sigma)&=&\{\sigma_{2,3}=(1,8)(2,4)(3,5)(6,7)\},\cr
C_{3\uar\dar}(\sigma)&=&C_{6\uar\dar}(\sigma)=\emptyset.
\end{array}$$

When  $\sigma_{i_s,i_t}$ is defined one has,
$$(R_{\sigma_{i_s,i_t}})_{i,j}=
\begin{cases}(R_{\sigma})_{i,j}-1, &{\rm if}\ i_t\leq i<i_s\ {\rm and}\ j_t\leq j<j_s;\cr
             (R_{\sigma})_{i,j},   & {\rm otherwise}.\cr
\end{cases}\eqno(***)$$
Thus, $\sigma_{i_s,i_t}<\sigma.$

In what follows we also need the following form of $(***)$:
$$(R_{\sigma_{i_s,i_t}})_{i,j}=
\begin{cases}(R_{\sigma_t^-})_{i,j}, &{\rm if}\ i>i_s\ {\rm
and}\ j<j_s;\cr
 (R_{\sigma})_{i,j},   & {\rm otherwise}.\cr
\end{cases}
 \eqno(***')$$
\begin{lemma}
Take $\sigma=(i_1,j_1)\ldots(i_k,j_k).$ Let $s\ :\ 1\leq s< k$ and
$t\ :\ t>s$ be such that $\sigma_{i_s,i_t}$ exists. Then
$\sigma_{i_s,i_t}\in D_1$.
\end{lemma}
\Pf
The proof is by induction on $n.$
 The minimal $n$ which permits such situation to occur is $n=4.$
In $\bS_4^2$ there  is a unique case of this type:
$\sigma=(1,4)(2,3)$ and $\sigma_{1,2}=(1,3)(2,4).$ It is obvious
that $\sigma_{1,2}\in D_1(\sigma).$ Assume the claim is true for
$n-1$.

Suppose $\sigma_{i_s,i_t}\in C_{i_s\uar\dar}(\sigma)$ for some $s.$
Let $\sigma\pr=(i_1\pr,j_1\pr)\ldots (i\pr_k,j\pr_k)$ be such that
$\sigma_{i_s,i_t}\leq \sigma\pr\leq \sigma.$ We must show that
either $\sigma\pr=\sigma$ or $\sigma\pr=\sigma_{i_s,i_t}.$

Note that $I(\sigma)=I(\sigma_{i_s,i_t})$ and
$J(\sigma)=J(\sigma_{i_s,i_t}).$ Thus, by lemma \ref{3.8b}
$I(\sigma\pr)=I(\sigma)$ and $J(\sigma\pr)=J(\sigma).$ We will
denote $I:=I(\sigma)$ and $J:=J(\sigma).$

Assume first that $1\not\in I(\sigma).$ Passing to $\pi_{2,n}$ we
get $\pi_{2,n}(\sigma)=\sigma,$\ $\pi_{2,n}(\sigma\pr)=\sigma\pr$\
and $\pi_{2,n}(\sigma_{i_s,i_t})=\sigma_{i_s,i_t}$ so that by the
induction hypothesis either $\sigma\pr=\sigma$ or
$\sigma\pr=\sigma_{i_s,i_t}.$

Now, let $1\in I(\sigma).$
\begin{itemize}
\item[(i)] Assume first that
$s>1.$ Then $\pi_{2,n}(\sigma_{i_s,i_t})=
(\pi_{2,n}(\sigma))_{i_s,i_t},$ thus, by induction hypothesis
$\pi_{2,n}(\sigma_{i_s,i_t})\in D_1((\pi_{2,n}(\sigma)))$ so that
either $\pi_{2,n}(\sigma\pr)=\pi_{2,n}(\sigma)$ or $
\pi_{2,n}(\sigma\pr)=(\pi_{2,n}(\sigma))_{i_s,i_t}.$ In both cases
$i_1\pr=1$ and $j_1\pr=j_1$ since $I(\sigma\pr)=I$ and
$J(\sigma\pr)=J.$ Thus, either $\sigma\pr=\sigma$ or
$\sigma\pr=\sigma_{i_s,i_t}.$
\item[(ii)] Now assume that $s=1.$ Consider $\pi_{2,n}(\sigma),$
$\pi_{2,n}(\sigma\pr)$ and $\pi_{2,n}(\sigma_{1,i_t}).$ Note that
$I(\pi_{2,n}(\sigma\pr))=I\setminus\{1\}=I(\pi_{2,n}(\sigma))=
I(\pi_{2,n}(\sigma_{1,i_t})).$

Further note that for any $m\ne t$ we get by the definition of
$M_{[1,j_s]}(\sigma)$ that either $j_m<j_t,$ or $i_m>i_t$ thus by
$(***)$\ $(R_{\sigma_{1,i_t}})_{p,r}=(R_{\sigma})_{p,r}$ for
$p=i_m,i_m+1$ and $r=j_m,j_m-1.$  This forces also
$(R_{\sigma\pr})_{p,r}=(R_{\sigma})_{p,r}$ for $p=i_m,i_m+1$ and
$r=j_m,j_m-1.$ By Lemma \ref{3.1}  $j\pr_m=j_m$ for any $m\ :\ m>1$
and $m\ne t.$ In other words for any $m\ne 1,t$ one has
$(i_m,j_m)\in\sigma\pr.$ Since $I(\sigma\pr)=I$ and $J(\sigma\pr)=J$
we get that either $(1,j_1),(i_t,j_t)\in\sigma\pr$ (that is
$\sigma\pr=\sigma$) or $(1,j_t),(i_t,j_1)\in\sigma\pr$ (that is
$\sigma\pr=\sigma_{1,i_t}$).
\end{itemize}
\QED
\subsection{}\label{3.8a}
\begin{prop}
 For $\sigma=(i_1,j_1)\ldots(i_k,j_k)$ in $\bS_n^2$
$$D_1(\sigma)=\coprod\limits_{s=1}^k\sigma_{i_s\uar}\cup\coprod
\limits_{s=1}^k \sigma_{j_s\rar}\cup \coprod\limits_{s=2}^k
C_{i_s\uar\rar}(\sigma)\cup \coprod\limits_{s=1}^k
C_{i_s\uar\dar}(\sigma)$$
\end{prop}
\Pf
 Let us denote $D\pr(\sigma):=\coprod\limits_{s=1}^k\sigma_{i_s\uar}\cup\coprod
\limits_{s=1}^k \sigma_{j_s\rar}$ and
$D\prpr(\sigma):=\coprod\limits_{s=2}^k C_{i_s\uar\rar}(\sigma)\cup
\coprod\limits_{s=1}^k C_{i_s\uar\dar}(\sigma).$
 By lemmas \ref{3.8}-\ref{3.8iii} it remains to show that for any
$\sigma\pr=(i_1\pr,j_1\pr)\ldots (i\pr_k,j\pr_k)$ such that
$\sigma\pr<\sigma$ there exist $\sigma\prpr\in D\pr(\sigma)\cup
D\prpr(\sigma)$ such that $\sigma\prpr\geq \sigma\pr.$ This is a
straightforward but rather long computation. It involves the
comparison of $I(\sigma)$ and $J(\sigma)$ with $I(\sigma\pr)$ and
$J(\sigma\pr).$ As in the proof of proposition \ref{3.7b} the
argument is broken into the six claims.

Since $l(\sigma\pr)=l(\sigma)$ one has that $I(\sigma\pr)\cup
J(\sigma\pr)\ne I(\sigma)\cup J(\sigma)$ implies the existence of
$q\in I(\sigma\pr)\cup J(\sigma\pr)$ such that $q\not\in
I(\sigma)\cup J(\sigma).$
\begin{specia}
Assume that there exists $q\in I(\sigma\pr)$ such that $q\not\in
I(\sigma)\cup J(\sigma).$ Then there exists $s\in
M(\pi_{q+1,n}(\sigma))$ such that $\sigma_{i_s\uar}\geq \sigma\pr.$
\end{specia}
\Pf
 Note that $(R_{\sigma\pr})_{q,n}=(R_{\sigma\pr})_{q+1,n}+1$ and
$(R_\sigma)_{q,n}=(R_\sigma)_{q+1,n}.$ Since
$(R_{\sigma\pr})_{q,n}\leq (R_\sigma)_{q,n}$  this implies
$(R_{\sigma\pr})_{q+1,n}<(R_{\sigma})_{q+1,n}$ so that by
Proposition \ref{3.7b} there exists $s\in M(\pi_{q+1,n}(\sigma))$
such that $(\pi_{q+1,n}(\sigma))_s^-\geq \pi_{q+1,n}(\sigma\pr).$
Consider $\sigma_{i_s\uar}.$  Since $s\in M(\pi_{q+1,n}(\sigma))$
and $q\not\in I(\sigma)$ one has $\sigma_{i_s\uar}\ne\emptyset.$
Moreover $m\geq q.$ Hence,
$$(R_{\sigma_{i_s\uar}})_{i,j}
\begin{cases} \geq (R_{(\pi_{q+1,n}(\sigma))_s^-})_{i,j}, & {\rm if}\ i\geq q+1;\cr
                                =(R_\sigma)_{i,j}, & {\rm if}\ i\leq q.\cr\end{cases}$$
Thus, $\sigma_{i_s\uar}\geq \sigma\pr.$
  \QED

By symmetry of the rank matrix around the anti-diagonal we get the
same result for the case $q\in J(\sigma\pr).$ We need the complete
formulation of this result  in what follows.
\begin{specia}
Assume that there exists $q\in J(\sigma\pr)$ such that $q\not\in
I(\sigma)\cup J(\sigma).$ Then there exists $r\in
M(\pi_{1,q-1}(\sigma))$ such that $\sigma_{j_r\rar}\geq \sigma\pr.$
\end{specia}

Now we have to consider the case $I(\sigma\pr)\cup J(\sigma\pr)=
I(\sigma)\cup J(\sigma).$ (Note that this does not imply
$I(\omega)\cup J(\omega)=I(\sigma)\cup J(\sigma)$ for all $\omega\
:\ \sigma\geq \omega\geq \sigma\prpr$). We show that in this case
there exists $\sigma\prpr\in D\pr(\sigma)\cup D\prpr(\sigma)$ such
that $\sigma\prpr\geq \sigma\pr.$  The proof is by induction. For
$n=4$ it holds trivially, as we have seen in \ref{3.8ii},
\ref{3.8iii}. Assume it for $n-1.$

\begin{specia}
 Assume that $I(\sigma\pr)\cup J(\sigma\pr)= I(\sigma)\cup J(\sigma).$
Let $n\not\in J(\sigma).$ Then there exists $\sigma\prpr\in
D\pr(\sigma)\cup D\prpr(\sigma)$ such that $\sigma\prpr\geq
\sigma\pr.$
\end{specia}
\Pf
Since $n\not\in J(\sigma), J(\sigma\pr)$  one has
$\pi_{1,n-1}(\sigma)=\sigma$ and $\pi_{1,n-1}(\sigma\pr)=\sigma\pr.$
Thus, by the induction hypothesis there exists $\sigma\prpr\in
D\pr(\pi_{1,n-1}(\sigma))\cup D\prpr(\pi_{1,n-1}(\sigma)$ such that
$\sigma\prpr\geq \pi_{1,n-1}(\sigma\pr).$ Considering $\sigma\prpr$
as an element of $\bS_n^2$ we get the result.
 \QED
Now assume that $n\in J(\sigma).$
Possibly changing the order of pairs we write
$\sigma=(i_1,j_1)\ldots (i_{k-1},j_{k-1})(i,n)$ and
$\sigma\pr=(i\pr_1,j\pr_1)\ldots(i\pr_{k-1},j\pr_{k-1})(i\pr,n).$
Let us start with the  case  $i\pr\ne i.$ In that case $i\in
I(\pi_{1,n-1}(\sigma\pr))\cup J(\pi_{1,n-1}(\sigma\pr))$ so that
either $i\in I(\pi_{1,n-1}(\sigma\pr))$  or $i\in
J(\pi_{1,n-1}(\sigma\pr)).$ We consider each case in a separate
claim.
 \begin{specia}
 Assume $I(\sigma\pr)\cup J(\sigma\pr)= I(\sigma)\cup J(\sigma).$
 Let $(i,n)\in \sigma,$ and $(i,n)\not\in\sigma\pr.$
Assume $i\in I(\sigma\pr).$  Then there exists $i_t\in I(\sigma)$
such that $i_t>i$ and $\sigma_{i,i_t}\geq \sigma\pr.$
\end{specia}
\Pf
 Since $i\not\in I(\pi_{1,n-1}(\sigma))$ one has
 $(R_{\sigma})_{i,n-1}=(R_{\sigma})_{i+1,n-1}.$
 Since $i\in I(\sigma\pr)$ and $(i,n)\not\in\sigma\pr$ one has $i\in I(\pi_{1,n-1}(\sigma\pr)).$
 Thus,
  $(R_{\sigma\pr})_{i,n-1}=(R_{\sigma\pr})_{i+1,n-1}+1.$
 Since $\sigma\pr<\sigma$, in particular, one has
 $(R_{\sigma\pr})_{i,n-1}\leq (R_{\sigma})_{i,n-1}.$
 Therefore,
 $(R_{\sigma\pr})_{i+1,n-1}<(R_{\sigma})_{i+1,n-1}.$  By
Proposition \ref{3.7b} there exists $t\in M(\pi_{i+1,n-1}(\sigma))$
such that $(\pi_{i+1,n-1}(\sigma))_t^-\geq
\pi_{i+1,n-1}(\sigma\pr).$ Note that
$(\pi_{i+1,n-1}(\sigma))_t^-=\pi_{i+1,n-1}(\sigma_t^-).$ In
particular, $(R_{\sigma\pr})_{p,q}\leq (R_{\sigma_t^-})_{p,q}$
whenever $p\geq i+1$ and $q\leq n-1.$

Note also that $t\in M_{[i,n]}(\sigma)$ so that $\sigma_{i,i_t}$
exists. By \ref{3.8iii} $(***')$
$$(R_{\sigma_{i,i_t}})_{p,q}=
\begin{cases}
 (R_{\sigma})_{p,q} &{\rm if}\ p\leq i\ {\rm or}\ q=n;\cr
 (R_{\sigma_t^-})_{p,q} & {\rm otherwise};\cr
 \end{cases}$$
Hence, $\sigma_{i,i_t}\geq \sigma\pr.$
 \QED
\begin{specia}
  Assume $I(\sigma\pr)\cup J(\sigma\pr)= I(\sigma)\cup J(\sigma).$
  Let $(i,n)\in\sigma$ and $i\in J(\sigma\pr)$. Then there exists $r\in
M(\pi_{1,i-1}(\sigma))$ such that $\sigma\prpr\geq\sigma\pr$ where
$\sigma\prpr=\sigma_{j_r\rar}$ or $\sigma\prpr=\sigma_{j_r,i}.$
\end{specia}
 \Pf

 Since $i\in J(\pi_{1,n-1}(\sigma\pr))$ and $i\not\in I(\pi_{1,n-1}(\sigma))\cup
 J(\pi_{1,n-1}(\sigma))$ by Claim 2
there exists $r\in M(\pi_{1,i-1}(\sigma))$ such that
$(\pi_{1,n-1}(\sigma))_{j_r\rar}\geq \pi_{1,n-1}(\sigma\pr).$

If there exists $r\in M(\pi_{1,i-1}(\sigma))$ such that $m<i$ then
$\sigma_{j_r\rar}\ne\emptyset$ and $\sigma_{j_r\rar}\geq \sigma\pr.$
Indeed,
$\pi_{1,n-1}(\sigma_{j_r\rar})=(\pi_{1,n-1}(\sigma))_{j_r\rar}$ so
that $(R_{\sigma_{j_r\rar}})_{p,q}\geq (R_{\sigma\pr})_{p,q}$
whenever $q<n$ and by $(*\pr)$ of \ref{3.8i}
$(R_{\sigma_{j_r\rar}})_{p,n}=(R_\sigma)_{p,n}\geq
(R_{\sigma\pr})_{p,n}$ for any $p.$

If for every $r\in M(\pi_{1,i-1}(\sigma))$ one has $m=i$  then
$\sigma_{j_r,i}$ exists for every $r.$ Let us take maximal such $r.$
Let us show that $\sigma_{j_r,i}\geq \sigma\pr.$
\begin{itemize}
\item[(i)]
$(R_{\sigma_{j_r,i}})_{p,q}=(R_{(\pi_{1,n-1}(\sigma))_{j_r\rar}})_{p,q}
\geq (R_{\sigma\pr})_{p,q}$ whenever $q<n.$

\item[(ii)] By \ref{3.8ii} $(**)$
$(R_{\sigma_{j_r,i}})_{p,n}=(R_\sigma)_{p,n}$ whenever either $p>i$
or $p\leq j_r.$ Thus, $(R_{\sigma_{j_r,i}})_{p,n}\geq
(R_{\sigma\pr})_{p,n}$ whenever either $p>i$ or $p\leq j_r.$

\item[(iii)] Since $i\in J(\sigma\pr), J(\sigma_{j_r,i})$ we get by the
form of the rank matrix
$(R_{\sigma\pr})_{i,n}=(R_{\sigma\pr})_{i+1,n}$ and
$(R_{\sigma_{j_r,i}})_{i,n}=(R_{\sigma_{j_r,i}})_{i+1,n}$ so that
$(R_{\sigma\pr})_{i,n}\leq (R_{\sigma_{j_r,i}})_{i,n}.$

\item[(iv)] Since $r$ is maximal and $m=i$ by definition of $\sigma_{j_r\rar}$
one has $t\in I(\sigma)$  for any $t\ :\ j_r<t<i,$ thus, also $t\in
I(\sigma_{j_r,i})$ for any such $t.$ Therefore, for any $t\ :\
j_r<t<i,$ one has
$(R_{\sigma_{j_r,i}})_{t,n}=(R_{\sigma_{j_r,i}})_{i,n}+(t-i).$ On
the other hand by the form of the rank matrix
$(R_{\sigma\pr})_{t,n}\leq (R_{\sigma\pr})_{i,n}+(t-i),$ so that
$(R_{\sigma\pr})_{t,n}\leq (R_{\sigma_{j_r,i}})_{t,n}$ for any $t\
:\ j_r<t<i.$
\end{itemize}
Together (i)-(iv) give us $\sigma_{j_r,i}\geq \sigma\pr.$
 \QED
\begin{specia}
Assume $I(\sigma\pr)\cup J(\sigma\pr)= I(\sigma)\cup J(\sigma).$
  Let $(i,n)\in\sigma,\sigma\pr.$
Then there exists $\sigma\prpr\in D\pr(\sigma)\cup D\prpr(\sigma)$
such that $\sigma\prpr\geq \sigma\pr.$
\end{specia}
\Pf
 Since $I(\pi_{1,n-1}(\sigma))\cup J(\pi_{1,n-1}(\sigma))=
I(\pi_{1,n-1}(\sigma\pr))\cup J(\pi_{1,n-1}(\sigma\pr))$ by the
induction hypothesis there exists $\hat\sigma\in
D\pr(\pi_{1,n-1}(\sigma))\cup D\prpr(\pi_{1,n-1}(\sigma))$ such that
$\hat\sigma\geq \pi_{1,n-1}(\sigma\pr).$

Assume $i\not\in I(\hat\sigma)\cup J(\hat\sigma).$ Then
$\hat\sigma(i,n)\in \bS_n^2$ and $\sigma\pr\leq
\hat\sigma(i,n)<\sigma.$ Moreover,
\begin{itemize}
\item[(i)] If $\hat\sigma=(\pi_{1,n-1}(\sigma))_{i_q,i_r}$ then either $i_r<i$ or
$i_q>i$ so that $\hat\sigma(i,n)=\sigma_{i_q,i_r}.$
\item[(ii)] If $\hat\sigma=(\pi_{1,n-1}(\sigma))_{j_r,i_q}$ then either $i_q<i$ or
$j_r>i$ so that $\hat\sigma(i,n)=\sigma_{j_r,i_q}.$
\item[(iii)] If $\hat\sigma=(\pi_{1,n-1}(\sigma))_{i_s\uar}$ and
either $i_s<i$ or $m>i$ then $\hat\sigma(i,n)=\sigma_{i_s\uar}.$
\item[(iv)] If $\hat\sigma=(\pi_{1,n-1}(\sigma))_{j_s\rar}$ and
either $j_s>i$ or $m<i$ then $\hat\sigma(i,n)=\sigma_{i_s\rar}.$
\end{itemize}
Setting  $\sigma\prpr=\hat\sigma(i,n)$ we get the result.

It is left to consider the two cases when $i\in I(\hat\sigma)\cup
J(\hat\sigma).$
\begin{itemize}
\item[(i)] If $\hat\sigma=(\pi_{1,n-1}(\sigma))_{i_s\uar}$ and $m=i$
let $\sigma\prpr=\sigma_{i,i_s}.$ It exists since $s\in
M_{[i,n]}(\sigma).$ Moreover, $\sigma_{i,i_s}\geq \sigma\pr$ since
$(R_{\sigma_{i,i_s}})_{p,q}=(R_{\hat\sigma})_{p,q}$ whenever $q<n$
and $(R_{\sigma_{i,i_s}})_{p,n}=(R_\sigma)_{p,n}$ for any $p.$

\item[(ii)] If $\hat\sigma=(\pi_{1,n-1}(\sigma))_{j_s\rar}$ and $m(j_s)=i$
we need a more subtle analysis. First note that in this case
$(R_\sigma)_{1,i}>(R_{\sigma\pr})_{1,i}.$ Indeed,
$(R_{\hat\sigma})_{1,i}=(R_{\hat\sigma})_{1,i-1}+1$ and
$(R_{\sigma\pr})_{1,i}=(R_{\sigma\pr})_{1,i-1}.$ Thus,
$(R_{\hat\sigma})_{1,i-1}\geq (R_{\sigma\pr})_{1,i-1}$ provides
$(R_{\sigma})_{1,i}\geq
(R_{\hat\sigma})_{1,i}>(R_{\sigma\pr})_{1,i}.$ Since
$(R_\sigma)_{1,n}=(R_{\sigma\pr})_{1,n}$ there exists $j>i$ such
that
 \begin{itemize}
\item[(a)] $p\in J(\sigma\pr)$ implies $p\in J(\sigma)$
 for any $p\ :\ i<p<j;$
\item[(b)] $j\in J(\sigma\pr)$ and $j\not\in J(\sigma).$
\end{itemize}
Since $I(\sigma)\cup J(\sigma)=I(\sigma\pr)\cup J(\sigma\pr)$
statements (a) and (b) are equivalent to
\begin{itemize}
\item[(a')] $p\in I(\sigma)$ implies $p\in I(\sigma\pr)$ for any $p\ :\ i<p<j;$
\item[(b')]  $j\in I(\sigma)$ and $j\not\in I(\sigma\pr).$
\end{itemize}

Thus by (b') $(R_\sigma)_{j,n-1}=(R_\sigma)_{j+1,n-1}+1$
 and $(R_{\sigma\pr})_{j,n-1}=(R_{\sigma\pr})_{j+1,n-1}$ so that
$(R_\sigma)_{j,n-1}>(R_{\sigma\pr})_{j,n-1}.$ By (a') this provides
$(R_{\sigma})_{i+1,n-1}>(R_{\sigma\pr})_{i+1,n-1}.$

Thus, by Proposition \ref{3.7b} there exists $t$ such that
$\pi_{i+1,n-1}(\sigma\pr)\leq (\pi_{i+1,n-1}(\sigma))_t^-$ which
implies by the definition that $\sigma_{i,i_t}$ exists and
$\sigma\pr\leq \sigma_{i,i_t}.$
\end{itemize}
\QED

This completes the proof of Proposition \ref{3.8a}.
 \QED

\subsection{}\label{3.10}
The two previous propositions involve in their proofs long and
tricky combinatorics. At least I have not succeeded to find shorter
and simpler proofs. On the other hand the proof of the following
lemma is more straightforward.
\begin{lemma} For any $\sigma\in \bS_n^2$ and for any $\sigma\pr\in D_1(\sigma)\cup D_2(\sigma)$ one has
$\Bscr_{\sigma\pr}\st\ov\Bscr_\sigma.$
\end{lemma}
\Pf Consider $\sigma=(i_1,j_1)\ldots(i_k,j_k)$ and $\sigma\pr\in
D_1(\sigma)\cup D_2(\sigma).$ Since the metric topology has more
open sets than the Zariski topology it is enough to show that
$N_{\sigma\pr}\in \ov\Bscr_{\sigma}$ for the metric topology on
${\rm End}\,\Co^n$. We construct a one parameter family of matrices
$\{A_k\}\in \bB$ such that $\lim\limits_{k\rar \infty}A_k N_\sigma
A_k^{-1}=N_{\sigma\pr}$ in the metric topology. Each $A_k$ will be
an elementary upper-triangular matrix or a product of  elementary
upper-triangular matrices. We construct them explicitly. Let
$E_i(a)$ and $E_{i+j}(a)$ be the invertible matrices defined as
follows
$$(E_i(a))_{s,t}=\begin{cases} a, & {\rm if}\ s=t=i;\cr
                               \delta_{s,t}& {\rm otherwise.}\cr
\end{cases}\quad
(E_{i+j}(a))_{s,t}=\begin{cases} a, & {\rm if}\ s=i,\ t=j;\cr
                        \delta_{s,t}& {\rm otherwise.}\cr
\end{cases}$$
One checks that $N_{\sigma\pr}=\lim\limits_{k\rar \infty}A_k
N_{\sigma}A_k^{-1}$ where $A_k$ are given as follows
\begin{itemize}
\item[(i)] For $\sigma\pr=\sigma_s^-$ take $A_k=E_{i_s}(\frac{1}{k});$
\item[(ii)] For $\sigma\pr=\sigma_{i_s}$ take $A_k=E_{i_s}(\frac{1}{k})E_{m_s+i_s}(1);$
\item[(iii)] For $\sigma\pr=\sigma_{j_s}$ take $A_k=E_{i_s}(\frac{1}{k})
E_{j_s+m_{j_s}}(-k);$
\item[(iv)] For $\sigma\pr=\sigma_{j_r,i_s}$ take $A_k=E_{i_r}(-\frac{1}{k})E_{i_s}(\frac{1}{k})
E_{j_r+i_s}(1);$
\item[(v)] For $\sigma\pr=\sigma_{i_s,i_t}$ take $A_k=E_{j_t}(k) E_{j_s}(-{\frac{1}{k}}
E_{i_s+i_t}(k) E_{j_t+j_s}(\frac{1}{k}).$
\end{itemize}
\QED
\subsection{}\label{3.11}
As a corollary of Lemma \ref{3.6}, Propositions \ref{3.7b} and
\ref{3.8a}
 and Lemma \ref{3.10} we get
\begin{cor}
\begin{itemize}\item[(i)] Given $\sigma,\sigma\pr\in \bS_n^2$ one has $\Bscr_{\sigma\pr}\st\ov\Bscr_\sigma$
iff $\Vscr_{\sigma\pr}\st \Vscr_\sigma.$ In particular $\ov
\Bscr_\sigma=\coprod\limits_{\sigma\pr\leq\sigma}\Bscr_{\sigma\pr}.$
\item[(ii)] $\ov \Bscr_\sigma=\Bscr_\sigma\coprod\bigcup\limits_{\sigma\pr\in D_1(\sigma)\cup D_2(\sigma)}
\ov \Bscr_{\sigma\pr}.$
\end{itemize}
\end{cor}
\Pf First of all note that by \ref{2.3} and $\bB$ stability of
$\Vscr_\sigma$ each $\Vscr_\sigma$ is a union of some of the
$\Bscr_{\sigma\pr}.$
\begin{itemize}
\item[(i)] By Lemma \ref{3.6} $\Bscr_{\sigma\pr}\st\ov\Bscr_\sigma$ implies
$\Vscr_{\sigma\pr}\st \Vscr_\sigma.$ On the other hand, if
$\Vscr_{\sigma\pr}\st \Vscr_\sigma$ then there exists a chain
$\sigma=\eta_1<\eta_2, \ldots <\eta_s=\sigma\pr$ such that
$\eta_{i+1}\in D_1(\eta_i)\cup D_2(\eta_i)$ for every $1\leq i\leq
s-1.$ Thus, by Lemma \ref{3.10} one has
$\Bscr_{\sigma\pr}\st\ov\Bscr_{\eta_{s-1}}\st\cdots\st\ov\Bscr_\sigma.$

\item[(ii)] By Lemma \ref{3.10}
$\bigcup\limits_{\sigma\pr\in D_1(\sigma)\cup
D_2(\sigma)}\ov\Bscr_{\sigma\pr}\st\ov\Bscr_\sigma.$ On the other
hand, if $\ov\Bscr_{\sigma\pr}\subsetneq\ov\Bscr_\sigma$ then
$\sigma\pr<\sigma$ just by Lemma \ref{3.6}. Thus, by Propositions
\ref{3.7b} and \ref{3.8a} there exists $\sigma\prpr\in
D_1(\sigma)\cup D_2(\sigma)$ such that $\sigma\pr\leq \sigma\prpr$.
Thus, by (i) $\Bscr_{\sigma\pr}\st\ov\Bscr_{\sigma\prpr}.$ This
concludes the proof.
\end{itemize}
\QED
\subsection{}\label{3.12}
Now we can prove theorem \ref{3.5}. Let us recall it
\begin{theorem} For any $\sigma\in\bS_n^2$ one has
$$\ov\Bscr_\sigma=\Vscr_\sigma=\coprod\limits_{\sigma\pr\leq\sigma}\Bscr_{\sigma\pr} .$$
\end{theorem}
\Pf By \ref{3.4} $\ov\Bscr_\sigma\st \Vscr_\sigma.$ Let us show that
$\dim \Vscr_\sigma=\dim \Bscr_\sigma$ and that $\Bscr_\sigma$ is the
only component in $\Vscr_\sigma$ of maximal dimension. Indeed, since
$\Bscr_\sigma\st \Vscr_\sigma$ we get that $\dim \Vscr_\sigma\geq
\dim \Bscr_\sigma.$ For any $\sigma\pr$ such that $\dim
\Bscr_{\sigma\pr}\geq \dim\Bscr_\sigma$ we have
$\Bscr_{\sigma\pr}\not\in \ov\Bscr_\sigma.$ Thus by \ref{3.11}(i)
one has $\sigma\pr\not<\sigma.$ That means that there exists $i,j$
such that $(R_{\sigma\pr})_{i,j}>(R_{\sigma})_{i,j},$ so that
$\Bscr_{\sigma\pr}\not\st \Vscr_\sigma.$ Hence, $\dim \Vscr_\sigma
=\dim \Bscr_\sigma$ and $\Bscr_\sigma$ is the only component of this
dimension in $\Vscr_\sigma.$ This together with \ref{3.7a} and
\ref{3.8} gives:
$$\Vscr_\sigma=\Bscr_\sigma\coprod\bigcup\limits_{\sigma\pr\in D_1(\sigma)
\cup D_2(\sigma)}\Vscr_{\sigma\pr}.$$ The result follows by
induction on order. Indeed, it is immediate that
$\Vscr_\sigma=\ov\Bscr_\sigma=\Bscr_\sigma\coprod \{0\}$ for
$\sigma=(1,n)$ which is the minimal non-zero orbit. Now assume that
this is true for all $\sigma\pr>\sigma$ then
 $$\Vscr_\sigma=\Bscr_\sigma\coprod\bigcup\limits_{\sigma\pr\in D_1(\sigma)\cup D_2(\sigma)}
\Vscr_{\sigma\pr}= \Bscr_\sigma\coprod \bigcup\limits_{\sigma\pr\in
D_1(\sigma)\cup D_2(\sigma)}\ov\Bscr_{\sigma\pr}
=\ov\Bscr_\sigma=\coprod\limits{\sigma\pr\geq\sigma}\Bscr_{\sigma\pr}.$$
\QED
\section {\bf Results on orbital varieties of nilpotent order 2}
\subsection{}\label{2.7}
Let us apply the results of section 3 to the study of orbital
varieties of nilpotent order 2 in $\gs\gl_n.$ First we recall some
facts about such orbital varieties and then in subsections
\ref{3.13}, \ref{3.14} formulate and prove the corollaries of
section 3 for orbital varieties.

Recall the notion of a Young tableau from \ref{1.4}. We  fill the
boxes of
 $D_\lam\in \bD_n$ with $n$ distinct positive integers (instead of integers $1,\ldots,n$)
in such a manner that entries increase in rows from left to right
and in columns from top to bottom. We call such an array a Young
tableau or simply a tableau. If the numbers in a Young tableau form
the set of integers from $1$ to $n,$ then the tableau is called a
standard Young tableau. Again, for any Young tableau $T$ associated
to $D_\lam$ we put $\sh(T)=\lam.$
\par
Let $\bT_n$ denote the set of all standard Young tableaux of size
$n.$
\par
Let $\bT_n^2\in \bT_n$ be the subset of standard Young tableaux with
two columns. Let $T=(T_1,T_2),$ where
$T_1=\left(\begin{array}{c}t_{1,1}\cr\vdots\cr
t_{l,1}\cr\end{array}\right)$ is the first column of $T$ and
$T_2=\left(\begin{array}{c}t_{1,2}\cr\vdots\cr
t_{k,2}\cr\end{array}\right)$ is the second column of $T.$
\par
Set $\sigma_{\sr T}=(i_1,j_1)\ldots (i_k,j_k)$ where $j_s=t_{s,2};\
i_1=t_{1,2}-1,$ and $i_s=\max\{d\in
T_1\setminus\{i_1,\ldots,i_{s-1}\}\ |\ d<j_s\}$ for any $s>1.$ Note
that now we order the cycles of $\sigma$ so that the second entries
increase. For example, take
$$T=\vcenter{
\halign{& \hfill#\hfill \tabskip4pt\cr \multispan{5}{\hrulefill}\cr
\ssa \vb & 1 & & 4 &\ts\vb\cr\vsa &&&&\cr \ssa \vb & 2 & & 5 &
\ts\vb\cr\vsa &&&&\cr \ssa \vb & 3 & & 7 & \ts\vb\cr\vsa &&&&\cr
\ssa \vb & 6 & & 8 & \ts\vb\cr\vsa \multispan{5}{\hrulefill}\cr}},$$
Then $\sigma_{\sr T}=(3,4)(2,5)(6,7)(1,8).$
\par
Let us denote $\Bscr_T:=\Bscr_{\sigma_T}.$ Recall the notion of
$\Uscr_T$ from \ref{1.4}. As we explained there if an orbital
variety is associated to a nilpotent orbit of nilpotent order 2,
then it contains a dense $\bB-$orbit. Explicitly by
\cite[4.13]{Msmith} one has
\begin{prop} For $T\in \bT_n^2$ one has $\ov\Bscr_T=\ov\Uscr_T.$
\end{prop}
\subsection{}\label{2.8}
Let us induce geometric and algebraic orders, defined in \ref{2.6}
from $\bS_n^2$ to $\bT_n^2.$ For $T,S\in \bT_n^2$ set $T\go S$ if
$\sigma_{\sr T}\go \sigma_{\sr S}$ and set $T\ao S$ if $\sigma_{\sr
T}\ao \sigma_{\sr S}.$ As a straightforward corollary of
\cite[3.7]{Mnew} we get that these two orders coincide, that is
\begin{theorem}
For $T,S\in \bT_n^2$ one has $T\go S$ iff $T\ao S.$
\end{theorem}
Since these two orders coincide on $\bT_n^2$ we can define a partial
order on $\bT_n^2$ by taking $T\leq S$ if $T\go S.$
\subsection{}\label{2.9}
For $T\in\bT_n^2$ set $\Dscr(T)$ to be  what is called in
combinatorics the cover of $T$, that is
$$\Dscr(T):=\{S\ |\ S>T\ {\rm and\ if\ } S\geq U> T\ {\rm then\ } U=S\}.$$
In \cite[3.18]{Mnew} $\Dscr(T)$ is described explicitly.
\par
Given $T=(T_1,T_2)$ where
$T_1=\left(\begin{array}{c}t_{1,1}\cr\vdots\cr
t_{l,1}\cr\end{array}\right)$ is the first column of $T$ and
$T_2=\left(\begin{array}{c}t_{1,2}\cr\vdots\cr
t_{k,2}\cr\end{array}\right)$ is the second column of $T.$ For
$t_{s,2}$ let $j_s:=\max\{j\ |\ t_{j,1}<t_{s,2}\}.$ Set $T\langle
t_{s,2}\rangle$ to be the tableau obtained from $T$ by moving
$t_{s,2}$ from the second column to the first one, that is $T\langle
t_{s,2}\rangle:=(T\pr_1,T\pr_2),$ where
$$T\pr_1=\left(\begin{array}{c}t_{1,1}\cr \vdots\cr t_{j_s,1}\cr t_{s,2}\cr t_{j_s+1,1}\cr\vdots
\cr t_{l,1}\cr\end{array}\right)\qquad {\rm and}\qquad
T\pr_2=\left(\begin{array}{c}t_{1,2}\cr\vdots\cr t_{s-1,2}\cr
t_{s+1,2}\cr \vdots\cr t_{k,2}\cr
\end{array}\right).$$
For example, take $T$ from \ref{2.7}. Then
$$T\langle 4\rangle=\vcenter{
\halign{& \hfill#\hfill \tabskip4pt\cr \multispan{5}{\hrulefill}\cr
\ssa \vb & 1 & & 5 &\ts\vb\cr\vsa &&&&\cr \ssa \vb & 2 & & 7 &
\ts\vb\cr\vsa &&&&\cr \ssa \vb & 3 & & 8 & \ts\vb\cr\vsa
&&\multispan{3}{\hrulefill}\cr \ssa \vb & 4 & \ts\vb\cr\vsa &&\cr
\ssa \vb & 6 & \ts\vb\cr\vsa \multispan{3}{\hrulefill}\cr}},\
T\langle 5\rangle=\vcenter{ \halign{& \hfill#\hfill \tabskip4pt\cr
\multispan{5}{\hrulefill}\cr \ssa \vb & 1 & & 4 &\ts\vb\cr\vsa
&&&&\cr \ssa \vb & 2 & & 7 & \ts\vb\cr\vsa &&&&\cr \ssa \vb & 3 & &
8 & \ts\vb\cr\vsa &&\multispan{3}{\hrulefill}\cr \ssa \vb & 5 &
\ts\vb\cr\vsa &&\cr \ssa \vb & 6 & \ts\vb\cr\vsa
\multispan{3}{\hrulefill}\cr}},\ T\langle 7\rangle =\vcenter{
\halign{& \hfill#\hfill \tabskip4pt\cr \multispan{5}{\hrulefill}\cr
\ssa \vb & 1 & & 4 &\ts\vb\cr\vsa &&&&\cr \ssa \vb & 2 & & 5 &
\ts\vb\cr\vsa &&&&\cr \ssa \vb & 3 & & 8 & \ts\vb\cr\vsa
&&\multispan{3}{\hrulefill}\cr \ssa \vb & 6 & \ts\vb\cr\vsa &&\cr
\ssa \vb & 7 & \ts\vb\cr\vsa \multispan{3}{\hrulefill}\cr}},\
 T\langle 8\rangle =\vcenter{
\halign{& \hfill#\hfill \tabskip4pt\cr \multispan{5}{\hrulefill}\cr
\ssa \vb & 1 & & 4 &\ts\vb\cr\vsa &&&&\cr \ssa \vb & 2 & & 5 &
\ts\vb\cr\vsa &&&&\cr \ssa \vb & 3 & & 7 & \ts\vb\cr\vsa
&&\multispan{3}{\hrulefill}\cr \ssa \vb & 6 & \ts\vb\cr\vsa &&\cr
\ssa \vb & 8 & \ts\vb\cr\vsa \multispan{3}{\hrulefill}\cr}}.$$ Note
that $T\langle t_{s,\sr 2}\rangle$ is always a Young tableau.
\par
Recall $\sigma_s^-$ from \ref{3.7a}.  Let $\sigma_{\sr
T}=(i_1,t_{1,2})\ldots (i_k, t_{k,2}).$ Let
$$\sigma_{\sr T\langle t_{s,2}\rangle} =(i\pr_1,t_{1,2})\ldots(i\pr_{s-1},t_{s,2})(i\pr_{s+1},t_{s+1,2})\ldots(i\pr_k, t_{k,2}).$$
Note that for any $p<s$ one has $i\pr_p=i_p.$ For $p>s$ one has
$i\pr_p=i_p$ if and only if $i\pr_p>t_{s,\sr 2}.$ In other words,
$\sigma_{\sr T\langle t_{s,2}\rangle} =(\sigma_{\sr T})_s^-$ if and
only if for any $p>s$ one has $i_p>t_{s,\sr 2}.$ This note permits
us to reformulate \cite[3.18]{Mnew} as follows
\begin{theorem} For any $T\in \bT_n^2$ one has
$$\Dscr(T)=\{T\langle t_{s,2}\rangle\ |\ \sigma_{\sr T\langle t_{s,2}\rangle} =
(\sigma_{\sr T})_s^- \}$$
\end{theorem}
For example, for $T$ from \ref{2.7} one has
$$\Dscr(T)=\{T\langle 8\rangle\}=\left\{\vcenter{
\halign{& \hfill#\hfill \tabskip4pt\cr \multispan{5}{\hrulefill}\cr
\ssa \vb & 1 & & 4 &\ts\vb\cr\vsa &&&&\cr \ssa \vb & 2 & & 5 &
\ts\vb\cr\vsa &&&&\cr \ssa \vb & 3 & & 7 & \ts\vb\cr\vsa
&&\multispan{3}{\hrulefill}\cr \ssa \vb & 6 & \ts\vb\cr\vsa &&\cr
\ssa \vb & 8 & \ts\vb\cr\vsa
\multispan{3}{\hrulefill}\cr}}\right\}$$
\subsection{}\label{3.13}
As a straightforward corollary of Proposition \ref{2.7} and Theorem
\ref{3.5} we get
\begin{cor} For any  $T\in\bT_n^2$ the ideal of definition $I(\Uscr_T)$
is the radical of $\Iscr_{\Bscr_T}.$
 $$I(\Uscr_T)=\sqrt{\Iscr_{\Bscr_T}}.$$
\end{cor}
\Pf Indeed, by \ref{2.7} $\ov\Uscr_T=\ov\Bscr_T$ and by \ref{3.5}
$\Vscr(\Iscr_{\sigma_T})=\ov\Bscr_T$ thus
$I(\Uscr_T)=\sqrt{\Iscr_{\Bscr_T}}.$ \QED
\subsection{}\label{3.14}
We get also
\begin{theorem} For $T\in\bT_n^2$ one has
$$\ov\Uscr_T=\bigcup\limits_{S\geq T}\Uscr_S$$
\end{theorem}
\Pf Since $\Bscr_T$ is dense in $\Uscr_T$ we get
$$\Uscr_T=\ov\Bscr_T\cap\Oscr_T=\coprod\limits_{{\sigma\pr\leq \sigma}\atop {l(\sigma\pr)=l(\sigma)}}
\Bscr_{\sigma\pr}.$$ On the other hand by theorem \ref{3.5}
$$\ov\Bscr_T=\coprod\limits_{\sigma\pr\leq \sigma}\Bscr_{\sigma\pr}=
\coprod\limits_{{\sigma\pr\leq \sigma}\atop
{l(\sigma\pr)=l(\sigma)}}
\Bscr_{\sigma\pr}\coprod\coprod\limits_{{\sigma\pr< \sigma}\atop
{l(\sigma\pr)<l(\sigma)}} \Bscr_{\sigma\pr}.$$ Further, by
definition of $D_2(\sigma)$ we get that for any $\sigma\pr$ such
that $\sigma\pr<\sigma$ and $l(\sigma\pr)<l(\sigma)$ there exists
$\sigma\prpr\in D_2(\sigma)$ such that $\sigma\prpr\geq \sigma\pr.$
In other words, by theorem \ref{3.5},
$\Bscr_{\sigma\pr}\st\ov\Bscr_{\sigma\prpr}.$ Thus, we get
$$\ov\Uscr_T=\ov\Bscr_T=\Uscr_T\coprod\bigcup\limits_{\sigma\pr\in D_2(\sigma_{\sr T})}\ov \Bscr_{\sigma\pr}.$$
Let us show that $D_2(\sigma_{\sr T})=\{\sigma_{\sr S}\}_{S\in
\Dscr(T)}.$ Indeed, take  $T=(T_1,T_2)$  where
$T_1=\left(\begin{array}{c}t_{1,1}\cr\vdots\cr
t_{l,1}\cr\end{array}\right)$  and
$T_2=\left(\begin{array}{c}t_{1,2}\cr\vdots\cr
t_{k,2}\cr\end{array}\right).$ Then $\sigma_{\sr
T}=(t_{i_1,1},t_{1,2})\ldots (t_{i_k,1},t_{k,2})$ where
$t_{i_1,1}=\max\{t_{s,1}\ :\ t_{s,1}<a_{1,2}\}$ and
$t_{i_j,1}=\max\{t_{s,1}\ :\ t_{s,1}\in
T_1\setminus\{t_{i_1,1},\ldots,t_{i_{j-1},1}\},\ t_{s,1}<t_{j,k}\}$
whenever $j>1.$ Since now the cycles of $\sigma_{\sr T}$  are
ordered so that the second entries increase, one has $j\in
M(\sigma_{\sr T})$ $\Longleftrightarrow$
$t_{i_j,1}=\max\{t_{i_1,1},\ldots,t_{i_j,1}\}$ so that $j\in
M(\sigma_{\sr T})$ if and only if  $(\sigma_{\sr
T})_j^-=\sigma_{T\langle t_{j,2}\rangle} $ for some $T\langle
t_{j,2}\rangle \in \Dscr(T).$ Thus, we get $\ov
\Bscr_T=\Uscr_T\coprod\bigcup\limits_{S\in \Dscr(T)}\ov\Uscr_S.$

The proof is completed by the induction on rank of $\Uscr_T.$ If
$\lam(T)=(n-1,1)$ that is $\Rank \Uscr_T=1$ it is obvious that
$\ov\Uscr_T=\Uscr_T\coprod\{0\}=\bigcup\limits_{S\geq T}\Uscr_S.$
Assume that this is true for $S$ of rank less than $k.$ Let $T$ be
of rank $k.$ Note that that for any $S>T$ one has $\Rank(S)\leq
k-1.$ Thus,
$$\ov\Uscr_T=\Uscr_T\coprod\bigcup\limits_{S\in \Dscr(T)}\ov\Uscr_S= \Uscr_T\coprod\
\bigcup\limits_{S\in \Dscr(T)}\left( \bigcup\limits_{U\geq
S}\Uscr_U\right)=\bigcup\limits_{S\geq T}\Uscr_S.$$ \QED
\bigskip
\parno
\centerline{ INDEX OF NOTATION}
\bigskip
\parno
\begin{tabular}{llll}
\ref{1.1} & $\bG,\ \bB,\ \bB_n,\ \gog,\ \gn,\ \gn_n,\ \gh,\ \gn^-,\
                                           \gn^-_n,$& \ref{2.5} & $\lambda\leq \mu$\cr

         & $S(*),\ \Oscr_u,\ \Bscr_u,\ \Xscr_2,\ \bS_n,\ \bS_n^2,\
                                       \Iscr_\Bscr,\ \Vscr(I)$&\ref{2.5a}& $x_{i,j}$\cr

\ref{1.1a}& $\Uscr$&\ref{2.6} & $\sigma\pr\go\sigma,\
                                                                 \sigma\pr\ao\sigma$\cr

\ref{1.3} & $\Wscr\go\Uscr$&\ref{2.11}& $(i,j)\in\sigma,\
                         I(\sigma),\ J(\sigma),\
                                                 S_{i,j}(\sigma),\ \pi_{i,j}(\sigma)$\cr

\ref{1.4} & $\lam^*,\ \lam,\ \Oscr_\lam,\ D_\lam,\ T,\ \sh(T),\
                                   D(u),$&\ref{3.1} & $R_u,\ R_{\Bscr_u},\ R_\sigma$\cr

          & $\pi_{1,i},\ D_i(u),\ \varphi(u),\ \nu_{\sr T},\ \Uscr_T$ &\ref{3.2} &
                                                                           $\bR_n^2$\cr

\ref{1.6} & $\pi_{i,n},\ D^i(u),\ \vartheta(u),\ \nu^{\sr T},$
             &\ref{3.4} & $I_2,\ \Iscr_{R_\sigma},\ \Iscr_\sigma,\
                                                     \Vscr_{R_\sigma},\ \Vscr_\sigma$\cr

          & $\pi_{i,j}(T),\ D_{i,j}(T),\ \nu\pr_{\sr T},\ \Iscr_T$&\ref{3.5} & $R\leq R\pr,\
                                                                \sigma\leq \sigma\pr$\cr

\ref{2.1} & $P(n),\ \bD_n$&\ref{3.7} & $D_1(\sigma),\ D_2(\sigma),\
                                                                         M(\sigma),$\cr

\ref{2.1a}& $\sh(u),\ \sh(\Oscr_u),\ \Oscr_{(l,k)}$&
                                                              \ref{3.7a}&$\sigma_s^-$\cr

\ref{2.2} & $e_{i,j},\ R,\ R^+,\ \Pi,\ \al_{i,j}, \al_i, X_{i,j},\
                                         \bP_{\al_i},$&\ref{3.8} & $\sigma_{i_s\uar}$\cr

     & $\bP_{\langle  i,j\rangle} ,\ M_{\langle  i,j\rangle} ,\ L_{\langle  i,j\rangle} ,\ \bB_{\langle  i,j\rangle} ,\
     \gm_{\sr \langle i,j\rangle} ,\qquad$&
                                                        \ref{3.8i}& $\sigma_{j_s\rar}$\cr

     & $\gl_{\sr \langle i,j\rangle} ,\ \gn_{\sr \langle i,j\rangle} ,\ \pi_{i,j},\ D_{i,j}(u),\ D_u,\ D_{\Bscr_u}$&
                                \ref{3.8ii}& $\sigma_{j_r,i_s},\ C_{i_s\uar\rar}(\sigma)$\cr

\ref{2.2a}& $W,\ s_{\al},\ s_i,\ S_{\langle  i,j\rangle} ,\
                          D_{\langle  i,j\rangle} ,$&\ref{3.8iii}& $\sigma_{i_s,i_t},\ C_{i_s\uar\dar}$\cr

          & $\pi_{i,j}:\bS_n\rar \bS_{\langle  i,j\rangle} $&\ref{2.7} & $\bT_n,\ \bT^2_n,\ T_1,\ T_2,\
                                                               \sigma_{\sr T},\ \Bscr_T$\cr

 \ref{2.3} & $N_\sigma,\ \Bscr_\sigma$&\ref{2.8} & $T\leq S$\cr
 \ref{2.3a}& $\Oscr_\sigma,\ l(\sigma)$&\ref{2.9} & $\Dscr(T),\ T\langle i\rangle$\cr

\end{tabular}

\end{document}